\documentclass[a4paper,11pt,twoside]{paper}
\usepackage[utf8]{inputenc}
\usepackage[russian,greek,english]{babel}
\usepackage{amsmath,amssymb,amsthm,mathrsfs,mathtools,graphicx}
\usepackage{hyperref}
\newtheoremstyle{sans}{\parskip}{\parskip}{\itshape}
                       {0pt}{\bfseries\sffamily}{.}{ }{}
\newtheoremstyle{sansplain}{\parskip}{\parskip}{}
                       {0pt}{\bfseries\sffamily}{.}{ }{}

\theoremstyle{sans}
\newtheorem{prop}{Proposition}[section]

\newtheorem{thm}[prop]{Theorem}
\newtheorem{lem}[prop]{Lemma}

\theoremstyle{sansplain}
\newtheorem{rem}[prop]{Remark}
\newtheorem{defin}[prop]{Definition}

\makeatletter
\@addtoreset{equation}{section}
\makeatother

\newcommand{\CC}{\mathbb{C}}
\newcommand\Pb{\mathbb{P}}
\newcommand\Sc{\mathcal{S}}
\newcommand\Zb{\mathbb{Z}}

\renewcommand{\geq}{\geqslant}
\renewcommand{\leq}{\leqslant}

\newcommand\1{\leavevmode\hbox{\rm \small1\kern-0.35em\normalsize1}}
\newcommand\ind[1]{\1_{\{#1\}}} 

\def\D{\mathcal{D}}
\def\E{\mathbb{E}}
\def\H{\mathcal{H}}
\def\HB{\mathbb{H}}
\def\L{\mathcal{L}}
\def\M{\mathcal{M}}
\def\S{{\bf S}}
\def\DD{\displaystyle}
\def\ab{{\alpha\beta}}
\def\ba{{\beta\alpha}}
\def\egaldef{\stackrel{\mbox{\tiny def}}{=}}

\newcommand\wt[1]{\widetilde{#1}}

\newcommand{\fproof}{\hfill $\blacksquare$ }
\newenvironment{ienumerate}{

\begin{enumerate}}{\end{enumerate}}

\newenvironment{myitemize}{

\begin{itemize}}{\end{itemize}}

\begin{document}
\title{{\Large  \textbf{Functional equations as an important analytic method in  stochastic modelling and in combinatorics}}}
\author{Guy Fayolle\thanks{INRIA Paris-Rocquencourt, Domaine de Voluceau, BP 105, 78153 Le Chesnay Cedex - France. \hspace{1cm} Email: {\tt Guy.Fayolle@inria.fr}}}   
\date{\today}
\maketitle
\begin{abstract}
\noindent Functional equations (FE) arise quite naturally in the analysis of stochastic systems of different kinds\,: queueing and telecommunication networks, random walks, enumeration of planar lattice walks, etc. Frequently, the object is to determine the probability generating function of some positive random vector in $\Zb_+^n$. Although the situation $n=1$ is more classical, we  quote an interesting non local functional equation which appeared in modelling a divide and conquer protocol for  a muti-access broadcast channel. As for $n=2$,  we outline the theory reducing these linear FEs to boundary value problems of Riemann-Hilbert-Carleman type, with closed form integral solutions. Typical queueing examples analyzed over the last $45$ years are sketched. Furthermore, it is also sometimes possible to determine the \emph{nature of the functions} (e.g., rational, algebraic, holonomic), as  illustrated in a combinatorial context, where asymptotics are briefly tackled. For general situations (e.g., big jumps, or $n\ge3$), only prospective comments are made, because then no concrete theory exists.
 \end{abstract}


AMS $2000$ Subject Classification: Primary 60G50; secondary 30F10, 30D05

\section{Introduction}

It is now almost undisputable that analytic methods became  ubiquitous in probability. Indeed, the greek word \foreignlanguage{greek}{ἀναλυτικός} refers essentially to \emph{``the ability to analyze''}, which is nothing else but the very nature of science\,! During the last century, the impressive development of information sciences (in particular computer and telecommunications networks) led to a need of system modelling, which in turn highlighted a number of new interesting (and sometimes fascinating) mathematical objects. In the sequel, we shall  focus on a particular class of these objects, namely \emph{functional equations (FE)}. Needless to emphasize that this subject, in the past, attracted famous mathematicians, among them d'Alembert Euler, Abel, Cauchy, Riemann...

The paper is organized as follows. Section~\ref{sec:cra} presents an interesting non local FE of one single variable, encountered in the analysis of a \emph{divide and conquer algorithm}. In Section~\ref{sec:2var}, we consider  FE coming from the analysis of the invariant measure of random walks in the quarter plane, and show in particular how the can be reduced to Riemann-Hilbert-Carleman boundary value problems. 
Section~\ref{sec:ex} sketches four examples, emanating from queueing network models  resolved over the last forty five-years. In a combinatorial  context, Section~\ref{sec:counting} summarizes results concerning the nature of the counting generating functions, when  the \emph{group} of the random walk is finite. Some questions related to asymptotics are also tackled. The concluding Section~\ref{sec:misc} gives prospective remarks for more general situations (arbitrary big jumps, $n\ge3$, etc), noting that currently no concrete global theory exists.

\section{A non local  functional equation of one complex variable for a collision resolution algorithm (CRA)} \label{sec:cra} A huge literature has been devoted to FEs when the unknown  function depends on a single variable. In this respect, the reader may see the seminal prominent book by Kuczma \cite{KUC}. In a probabilistic context, we present a simple FE encountered in the analysis of a variety of the Capetanakis-Tsybakov-Mikhailov CRA, which is a \emph{divide and conquer} algorithm. All proofs can be found in \cite{FFH}.
\subsection{Specification of the CRA with continuous input}
\begin{enumerate}
\item  A single error-free channel is shared among many users which transmit messages of constant length (packets). Time is slotted and may be considered discrete. Users are synchronized with respect to the slots, and packets are transmitted at the beginning of slots only. Each slot is equal to the time required to transmit a packet (see the famous ALOHA network concept).
\item Each transmission is receivable by every user. Thus, when two or more
users transmit simultaneously, packets are said to \emph{collide} (interfere)
and none is received correctly: these collisions are treated as transmission
errors and each user must strive to retransmit its colliding packet
until it is correctly received. The users all employ the same algorithm for
this purpose, and have to resolve the contention without the benefit of
any other source of information on other users' activity save the
common channel.
\item Each user monitoring the channel knows, by the end of the slot, if that
slot produced a collision or not.
\item Each active user maintains a conceptual stack. At each slot end, he determines his position in the stack according to the following procedure (identical to all users, who are unable, however to communicate their stack state):
\begin{itemize}
\item[-] When an inactive user becomes active, it enters level $0$ in the stack. He will transmit at the nearest slot, and will always do so when at stack level $0$.
\item[-] After a non-collision slot, a user in stack level $0$ (there can be at
most one such user) becomes inactive, and all users decrease their stack level by $1$.
\item[-] After a collision slot, all users at stack level $i, i\ge1$ change to
level $i+1$. The users at level $0$ are split into two groups; one group remains at
level $0$, while the members of the other push themselves into level $1$. 
 This partition can be made on the basis of a Bernouilli trial, each user flipping a two-sided coin (independently of the other active users) : with probability $p$, he remains at level $0$, and with probability $q=1-p$ he pushes himself into level $1$.
  \end{itemize}
\item The numbers of new packets generated in each slot (i.e. the number of
new active users) form a sequence of i.i.d. random variables, denoted
by $X_i,i\ge1$, which follow a Poisson distribution with parameter $\lambda$.
\end{enumerate}
The \emph{collision resolution interval} (CRI), denoted in the sequal by $L_n$, is the time it takes to dispose of a group of $n$ colliders initially at level $0$.

\subsection{Functional equation for the generating function of the mean CRI}
The random variables $L_n$  satisfy the recursive relationship
\begin{equation}
\begin{cases} 
L_0=L_1=1,\\[0.2cm]
L_n = 1 + L_{I+X} + L_{n-I+Y}, \quad n \ge2,
\end{cases}
\end{equation}
where 
\begin{itemize}
\item $I$, the number of messages immediately retransmitted, follows the
binomial distribution $B(n, p)$;
\item $X$ is the number of new arrivals in that collision slot; 
\item $Y$ is the number of new arrivals in the slot following $L_{I+X}$. 
\end{itemize}
Moreover, $I,X,Y$ are supposed to be independent random variables 

Letting $\alpha_n\egaldef \E(L_n)$ and introducing 
\begin{equation}
\begin{cases} 
\DD \alpha(z)\egaldef \sum_{n\ge0}\alpha_n \frac{z^n}{n!}, \\[0.4cm]
\DD\psi(z)\egaldef e^{-z}\alpha(z),
\end{cases}
\end{equation}
we obtain the \emph{non local} FE
\begin{equation}\label{eq:FE}
\boxed{\psi(z) - \psi(\lambda + pz) - \psi(\lambda + qz)  = 1 - 2\psi(\lambda) e^{-z}(1 + Kz)},
\end{equation}
where 
\[
K= \frac{\exp(-\lambda p) - \exp(-\lambda q)}{\frac{\lambda}{q}\exp(-\lambda q) -\frac{\lambda}{p} \exp(-\lambda p)}.
\]

From now on, $p\ge q$ with $(p+q=1)$ and
$\sigma_1(z)\egaldef\alpha+pz$, $\sigma_2(z)\egaldef\alpha+qz$.

To solve \eqref{eq:FE}, we need to introduce a \emph{non-commutative iteration semigroup} $H$  of linear substitutions generated by $\sigma_1,\sigma_2$, where the semigroup operation is the composition of functions. The
identity of $H$ is denoted by $\varepsilon$, so that  $\varepsilon(z) = z, \forall z \in\CC$
(the complex plane) and any $\sigma\in H$ can be written in the form
\[
\sigma = \sigma_{i_1} \sigma_{i_2}\dots \sigma_{i_n}, \quad n\ge0, \ i_j\in\{1,2\}.
\]
Setting 
\[
 |\sigma|_1= \mathrm{card}\{j | i_j=1\}, \quad |\sigma|_2= \mathrm{card}\{j | i_j=2\}, \quad
 |\sigma|=  |\sigma|_1| + |\sigma|_2|,
 \]
we introduce the notation, valid for arbitrary complex numbers $\alpha,\beta$,
\[
(\alpha;\beta)^\sigma = \alpha^{|\sigma|_1|} \beta^{|\sigma|_2|}.
\]
By linearity, we have $\sigma(z)= \sigma(0) +(p;q)^\sigma z$.

\begin{figure}[!h]
\vspace{-1cm}
\includegraphics[width=\linewidth]{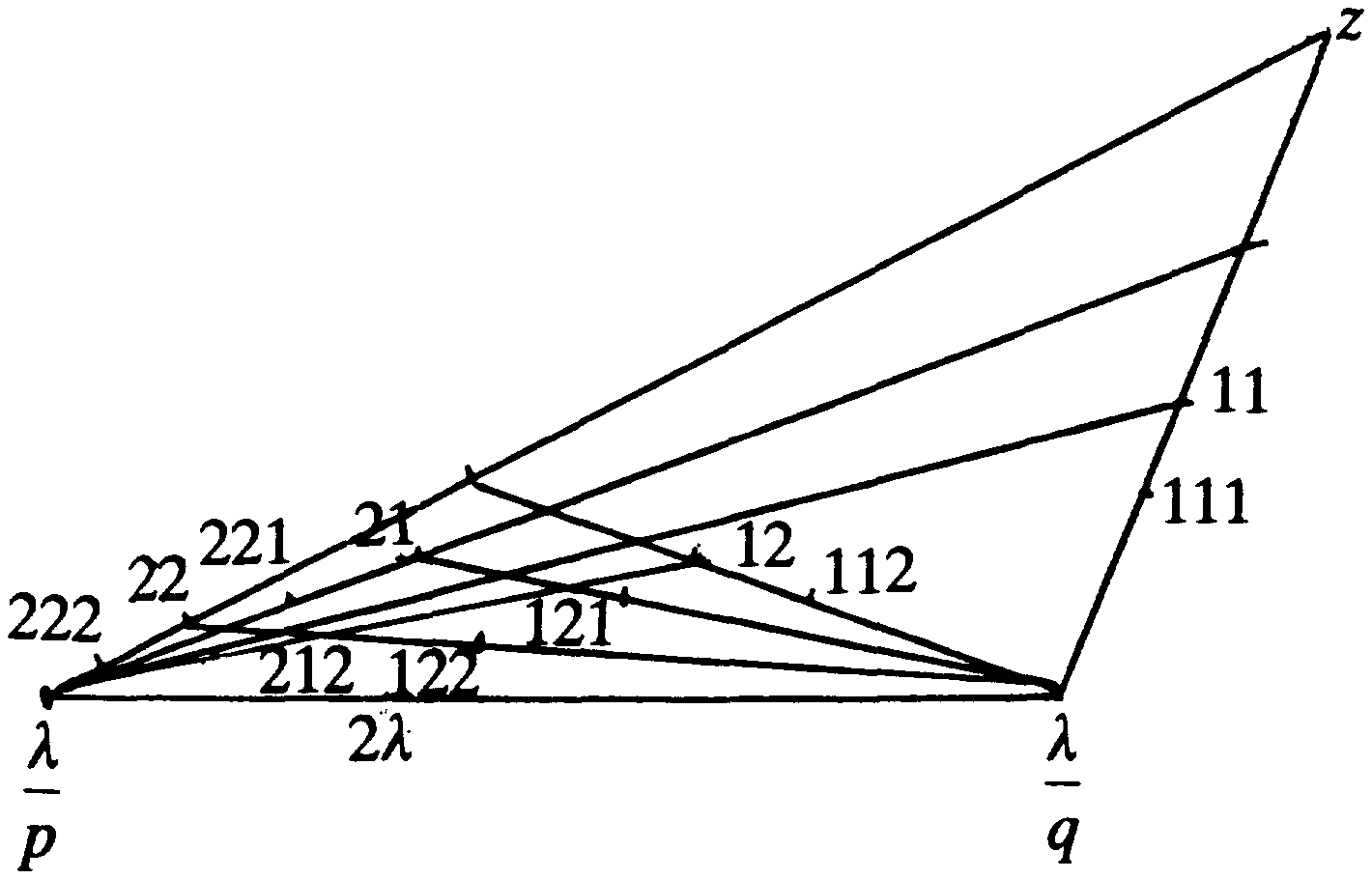}
\vspace{-13cm}
\caption{Successive transforms of a point $z$, in the case $p =2/3, q =1/3$.}
\end{figure}
Letting
\begin{equation*}
\begin{cases}
\DD g_n = (-1)^n \sum_\sigma\exp(-\sigma(0))(p^n;q^n)^\sigma, \\[0.4cm]
\DD k_n = (-1)^n \sum_\sigma\sigma(0)\exp(-\sigma(0))(p^n;q^n)^\sigma, \\[0.4cm]
\DD D(\lambda) = \sum_{n\ge2}[1-Kn)g_n +Kk_n]\frac{\lambda^n}{n!},
\end{cases}
\end{equation*}
we summarize the main results obtained in \cite{FFH} .
\medskip
Let 
\[
\mathbf{S}[f(\cdot);z]\egaldef \sum_{\sigma\in H} \big[f(\sigma(z))-f(\sigma(0)) -(p;q)^\sigma zf'(\sigma(0))\bigr].
\]
\begin{thm} \label{eq:FE1}
\[
\psi(z)= 1-\frac{2\, \mathbf{S}[e^{-u}(1+Ku);z]}{1+2D(\lambda)}.
\]
\end{thm}

\begin{thm} [Asymptotics]\label{eq:FFH}
The mean time to resolve $n$ collisions satisfies
\begin{equation}\label{eq:CRIn}
\alpha_n = \frac{2A}{1+2D(\lambda)} n +
 \frac{1}{1+2D(\lambda)} \sum_\chi a(\chi)n^{-\chi} + \mathcal{O}(n^{1-\eta}),
\end{equation}
for any sufficiently small $\eta> 0$, where the summation is extended to the $\chi$'s satisfying
\[
1-p^{-\chi}  - q^{-\chi}=0,  \  -1\leq \Re(\chi) < -1+\eta, \ \chi\ne -1.
\]
 The sum in the expression is a bounded fluctuating function, with an amplitude small compared to $nA$, typically less by several orders of magnitude, and the following properties hold.
 \begin{itemize}
\item $A$ is a complicated constant involving a Riemann-Stieltjes integral with respect to a measure having a nowhere differentiable  density.
 \item If \ $\DD \frac{\log p}{\log q}= \frac{d}{r}$ is \emph{rational}, i.e. with $\gcd (d,r)=1$, then 
 \begin{equation*}
\alpha_n = \frac{2A}{1+2D(\lambda)} n + n P(r\log_pn) + o(n^{1-\eta}),
\end{equation*}
with $P(u)$ a Fourier series of $u$ with mean value $0$. In this case  
$\DD\lim_{n\to\infty}\alpha_n/n$ does not exist.
\item If \ $\DD \frac{\log p}{\log q}$ is \emph{not rational}, then the sum in 
\eqref{eq:CRIn} is $o(n)$ and $\DD\lim_{n\to\infty}\alpha_n/n$ exists.
\end{itemize}
\end{thm}
\medskip

The main ingredients in the proof of Theorem~\ref{eq:FFH} are the exponential approximation [i.e. replace~$(1-a)^n$ by~$\exp(-an)$], together with a skillful use of Mellin's transforms, yielding the intermediate fundamental proposition.

 For $r(u)$ is any continuously differentiable function on $[0, \lambda/q]$,
define the  Dirichlet series 
\[
\omega(s) = \sum_{\sigma\in H} r(\sigma(0))(p^s;q^s)^\sigma.
\]
\begin{prop}
\begin{eqnarray*}
\omega(s)&=&\frac{1}{(s-1)}\frac{1}{h(p,q)(\frac{\lambda}{q}-\frac{\lambda}{p})}
 \int_{\lambda/p}^{\lambda/q}r(u)du \\
 &+&  \int_0^{\lambda/q}r'(u)w(u)du 
 + \frac{p\log^2p +q\log^2q}{h(p,q)(\frac{\lambda}{q}-\frac{\lambda}{p})} +o(s-1),
\end{eqnarray*}
where $h(p, q)\egaldef p\log p^{-1} + q\log q^{-1}$ is the entropy function and $w(u)$ is nowhere differentiable.
\end{prop}  
\begin{thm} [Ergodicity]
A necessary and sufficient condition to have a stable channel, i.e. $\alpha_n < \infty, \forall n$ finite, is $\lambda<\lambda_{max}$, where $\lambda_{max}$
is the first positive root of $1 +2D (\lambda) = 0$. The proof relies on standard results on Markov chains, using the stochastic interpretation of
$\psi(\lambda) = [1+2D(\lambda)^{-1}$. When $p = q = 1/2, \lambda_{max} = 0.3601$.
\end{thm}

\section{Functional equations of  two complex variables} \label{sec:2var}
In a probabilistic framework,  we consider a piecewise homogeneous random walk with sample paths in $\Zb_+^2$, the lattice in the positive quarter plane. In the strict interior of $\Zb_+^2$,  the size of the jumps is $1$,  and  $\{p_{ij}, |i|,|j| \leq 1\}$ will denote the generator of the process for this region. Thus a transition $(m,n)\to(m+i,n+j), m,n>0,$ can take place with probability $p_{ij}$, and
\[
\sum_{ |i|,|j| \leq 1} p_{ij} =1.
\]
No strong assumption is made about the boundedness of the upward jumps on the axes, neither at $(0,0)$. In addition, the downward jumps on the $x$ [resp.\ $y$] axis are bounded by $L$ [resp.\ $M$], where $L$ and $M$ are arbitrary finite integers. The basic problem is to determine the invariant measure $\{\pi_{i,j}, i,j \ge 0\}$, the  generating function of which satisfies the fundamental FE 
\begin{equation} \label{eq:eqfonc}
\boxed{Q(x,y) \pi(x,y) = q(x,y) \pi(x) + \wt{q}(x,y)\wt{\pi}(y) + \pi_0(x,y)}, 
\end{equation}
where $x,y$ belong to the complex plane $\CC$ with $|x|<1,|y|<1\}$, and
\begin{equation*}  
\label{A2-FE}
\begin{cases}
\pi(x,y) = \DD \sum_{i,j \geq 1} \pi_{ij} x^{i-1} y^{j-1}, \\[0.5cm]  
\pi(x) = \DD \sum_{i\geq L} \pi_{i0} x^{i-L}, \quad \wt{\pi}(y) = \sum_{j\geq M} \pi_{0j} y^{j-M}, \\[0.5cm] 
Q(x,y) = \DD xy \Bigg[ 1 - \sum_{i,j\in\Sc} p_{ij} x^i y^j \Bigg], \quad  \sum_{i,j\in\Sc} p_{ij} =1, \\[0.5cm] 
q(x,y) = \DD  x^L \Bigg[\sum_{i \geq -L, j \geq 0}
p'_{ij} x^i y^j - 1 \Bigg] \equiv x^L(P_{L0}(x,y)-1), \\[0.5cm]
\wt{q}(x,y) = \DD y^M \Bigg[\sum_{i \geq 0, j \geq -M} p''_{ij} x^i y^j - 1 \Bigg] 
\equiv y^M(P_{0M}(x,y)-1), \\[0.6cm]
 \pi_0(x,y) = \DD  \sum_{i=1} ^{ L-1}\pi_{i0} x^i\big[P_{i0}(x,y)-1\big] 
 +  \sum_{j=1} ^{ M-1}\pi_{0j} y^j\big[P_{0j}(x,y)-1\big] + \pi_{00}(P_{00}(xy)-1).
  \end{cases}
\end{equation*}
In equation \eqref{eq:eqfonc}, $\mathcal{S}$ is the set of allowed jumps, the unknown functions 
$\pi(x,y), \pi(x), \wt{\pi}(y)$ are sought to be \emph{analytic} in the region $\{(x,y)\in \mathbb{C}^{2} : |x|<1,|y|<1\}$, and \emph{continuous} on their respective boundaries. In addition, $q, \wt{q}, q_0, P_{i0}, P_{0j}, $ are given probability generating functions supposed to have suitable analytic continuations (as a rule, they are polynomials when the jumps are bounded). 

The polynomial  $Q(x,y)$ is often referred to as the \emph{kernel} of (\ref{eq:eqfonc}). 

Completely new approaches toward the solution of  the  problem
were discovered by the authors of the book \cite{FIM2}, the goal going far
beyond the mere obtention of an index theory for the quarter plane.
The main results can be summarized as follows.
\begin{enumerate}
\item  The first step, which is quite similar to a Wiener--Hopf factorization,
consists in considering the above equation on the algebraic curve
$\{Q(x,y)=0\}$ (which is {\em elliptic} in the generic situation), so that we
are then left  with an equation for two unknown functions of one
variable on this curve.
\item Next a crucial idea is to use \emph{Galois automorphisms} on this
algebraic curve. 
Let $\CC\,(x,y)$ be the field of rational functions in $(x,y)$ over $\CC$. Since $Q$ is assumed to be irreducible in the general case, the quotient field $\CC\,(x,y)$ denoted by 
 $\CC_Q(x,y)$ is also a field.
 \begin{defin}\label{def:gr}
The \emph{group of the random walk} is the Galois group  $\H=\langle \xi,\eta\rangle$ of automorphisms of   $\CC_Q(x,y)$ generated by  $\xi$ and $\eta$ given by
     \begin{equation*}
          \xi(x,y)=\Bigg(x,\frac{1}{y}\frac{\sum_{i} p_{i,-1}x^{i}}
          {\sum_{i} p_{i,1}x^{i}}\Bigg),\ \ \ \ \
          \eta(x,y)=\Bigg(\frac{1}{x}\frac{\sum_{j} p_{-1,j}y^{j}}
          {\sum_{j} p_{1,j}y^{j}},y\Bigg) .
     \end{equation*}
Here $\xi$ and $\eta$ are involutions satisfying  $\xi^2=\eta^2=I$.  Let
$\delta \egaldef\eta\xi$
denote their product, which is non-commutative except for $\delta^2=I$.Then $\H$ has a normal cyclic subgroup $\H_0=\{\delta^i, i\in\Zb\,\}$, which is finite or infinite, and 
$\H / \H_0$ is a group of order~$2$. Hence the group $\H$ is \emph{finite of order $2n$} if, and only if, 
\begin{equation}\label{eq:gr2n}
 \delta^n = I.
\end{equation}
\end{defin}
More information is obtained by using the fact
that the unknown functions $\pi$ and $\widetilde{\pi}$ depend
 solely on $x$ and $y$ respectively, i.e. they are invariant with
respect to $\xi$ and $\eta$ correspondingly.  It is then possible to prove that  $\pi$ and
$\widetilde{\pi}$ can be \emph{lifted} as meromorphic functions onto
the \emph{universal covering} of some Riemann surface $\S$. Here $\S$
corresponds to the algebraic curve $\{Q(x,y)=0\}$.  When \mbox{$g\egaldef$ \emph{the genus} of $\S$} is $1$ (resp. $0$), the universal covering is  the \emph{complex plane} $\CC$ (resp. the \emph{Riemann sphere}).
\item Lifted onto the universal covering, $\pi$ (and also
 $\widetilde{\pi}$) satisfies a system of non-local equations having
 the simple form
\[
\begin{cases}
\pi (t + \omega_1) = \pi (t), \quad \forall t \in \CC\,, \\ \pi (t +
\omega_3) = a(t)\pi (t) + b(t), \quad \forall t \in \CC\,,
\end{cases} \] 
where $\omega_1$ [resp. $\omega_3$] is a complex [resp. real]
constant. The solution can be presented in terms of infinite series
equivalent to {\em Abelian integrals}. The backward transformation
(projection) from the universal covering onto the initial coordinates
can be given in terms of uniformization functions, which, for $g=1$, are elliptic functions.
\item  Another direct approach to solving the fundamental equation
consists in working solely in the complex plane $\CC$. After making the
analytic continuation, it appears that the determination of $\pi$
reduces to a \emph{boundary value problem} (BVP), belonging to the Riemann--Hilbert--Carleman class, the basic form of which can be formulated as follows.

\medskip  Let $\mathscr{G}(\L)$ denote the interior of the domain bounded by
a simple smooth closed contour $\L$.

 \emph{Find a function $\Phi^+$ holomorphic in $\mathscr{G}(\L)$, the limiting
values of which are continuous  on the contour and satisfy the
relation}
\begin{equation}\label{eq:BVP}
 \Phi^+(\alpha(t)) = G(t) \Phi^+(t) + g(t), \quad t \in \L,
 \end{equation}
where
\begin{myitemize}
\item $g$, $G \in \HB_{\,\mu}(\L)$ ({\em H\"{o}lder condition} with
parameter $\mu$ on $\L$);
\item $\alpha$, referred to as a {\em shift} in the sequel, is a
function establishing a one-to-one mapping of the contour $\L$ onto
itself, such that the direction of traversing $\L$ is changed and
$$\alpha'(t) = \frac{d\alpha(t)}{dt} \in \HB_{\,\mu}(\L), \quad
\alpha'(t) \neq 0, \quad \forall t \in \L.$$
\end{myitemize}
In addition,  the function $\alpha$ is most frequently subject to the
 so-called \emph{Carleman condition}
$$\alpha(\alpha(t)) = t, \quad \forall t \in \L, \quad \text{where
 typically} \quad \alpha(t)=\overline{t}.$$
The advantage of this method resides in the fact that  solutions are given in terms of explicit integral-forms.
\item  Analytic continuation gives a clear understanding of possible
singularities and thus allows to derive the asymptotics of the functions.
\end{enumerate}
All these techniques work quite similarly for Toeplitz operators, and for other questions related  to random walks as well: transient behavior, first hitting time problem
\cite{GROM73}, calculating the Martin boundary, non spatially homogeneous walks, etc. 

\medskip For the sake of historical reference, it is worth quoting the pioneering work of V.A.~Malyshev  relating to points $1,2,3$,  which  was mainly settled
in the period 1968--1972 (see e.g.,  \cite{MALY70,MALY71c,MALY72a}).

\medskip  The method concerning points $1$ and $4$ was proposed in the
seminal study  \cite{FAIA79} carried out in 1976--1979, which was
widely referred to and followed up  in many other papers, until today. 
The three authors joined their efforts in the book \cite{FIM2}, where the reader can find a fairly comprehensive bibliography. 

\subsection{Summary of some general results (see \cite{FIM2}) }The multi-valued algebraic function $Y(x)$ solution of the polynomial equation 
\[
Q(x,Y(x))=0, \quad x\in\CC,
\]
is defined in the $\CC_x$-plane and has two branches $Y_0(x),Y_1(x)$. 
Rewrite for a while $Q(x,y)$ in the form
     \begin{equation}
     \label{eq:kernel}
          Q(x,y) = a(x) y^{2}+ b(x) y + c(x). 
     \end{equation}
\subsubsection{Genus $1$}\label{sec:genus1} In this case,  it  can be shown that $Y(x)$ has \emph{$4$ real branch points}, which are the roots of the discriminant $D(x)=b^2(x)-4a(x)c(x)$, two of them $x_1,x_2$ being located inside the unit disc $\D$.  In addition, there exists a \emph{uniformization} in terms of the Weierstrass $\wp$ function with periods $\omega_1, \omega_2$  depending on the parameters $p_{ij}$.

Clearly, exchanging $x$ and $y$, similar properties hold for the function $X(y)$ defined by $Q(X(y),y)=0$.

Let $[\underleftarrow{\overrightarrow{x_1x_2}}]$  stand for the
{\em contour} $[x_1 x_2]$, traversed from $x_1$ to $x_2$ along the
upper edge of the slit $[x_1 x_2]$ and then back to $x_1$, along the
lower edge of the slit. Similarly, $[\underrightarrow{\overleftarrow {x_1x_2}}]$
is defined by exchanging ``upper'' and ``lower''. Noting that on their respective cuts 
$Y_0(x)=\overline{Y}_1(x)$ and $X_0(y)=\overline{X}_1(y)$, one can set
 \begin{equation*}
\begin{cases} 
\L &= \ Y_0[\underleftarrow{\overrightarrow{x_1x_2}}] =
\overline{Y}_1[\underrightarrow{\overleftarrow{x_1x_2}}], \\ \L_{ext}
&= \ Y_0[\underleftarrow{\overrightarrow{x_3x_4}}] =
\overline{Y}_1[\underrightarrow{\overleftarrow{x_3x_4}}], \\ \M &= \
X_0[\underleftarrow{\overrightarrow{y_1y_2}}] =
\overline{X}_1[\underrightarrow{\overleftarrow{y_1y_2}}], \\ \M_{ext}
&= \ X_0[\underleftarrow{\overrightarrow{y_3y_4}}] = \overline{X}_1
[\underrightarrow{\overleftarrow{y_3y_4}}]. \end{cases}
\end{equation*}

\begin{figure}[!h]
\begin{center} 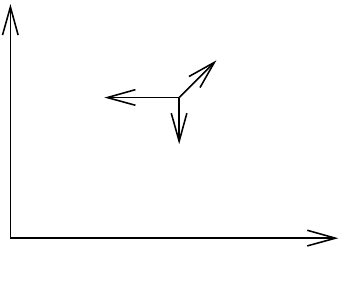 \end{center}
\begin{center} 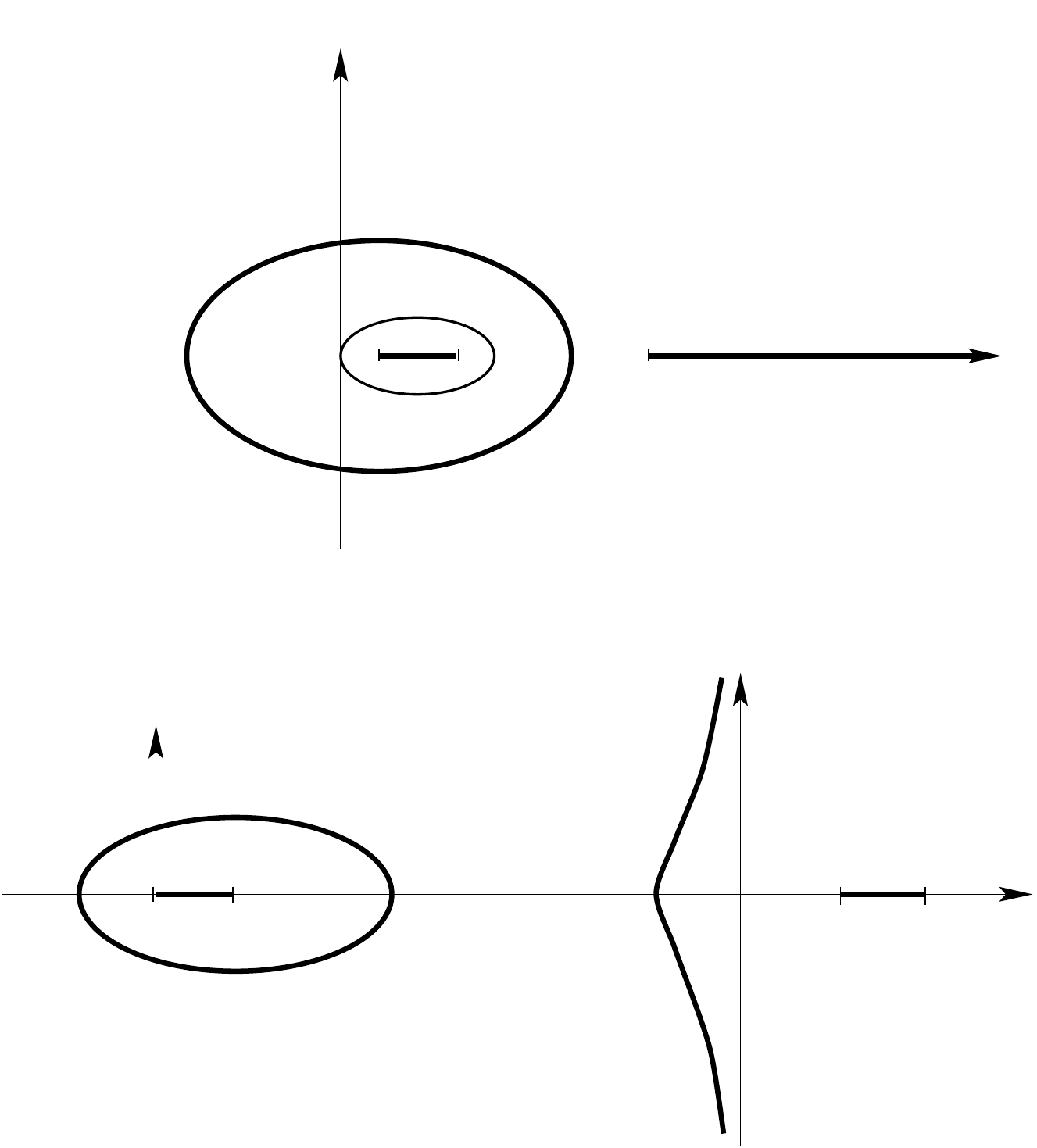 \end{center}
\caption{Example of the mappings of the cuts}  \label{fig:cuts}
\end{figure}

\begin{thm} [part of Theorem 5.3.3 in \cite{FIM2}]\mbox{ }
\begin{ienumerate}
\item The curves $\L$ and $\L_{ext}$ (resp. $\M$ and $\M_{ext}$)
are simple, closed and symmetrical about the real axis in the $\CC_y$
[resp. $\CC_x$] plane. They do not intersect if the group of the random
walk is not of order 4. When this group is of order 4, $\L$ and
$\L_{ext}$ [resp. $\M$ and $\M_{ext}$] coincide and form a circle
possibly degenerating into a straight line.  In the general case they
build the two components (possibly identical, in which case the circle must
be counted twice) of a quartic curve (see an example in figure~\ref{fig:cuts}).
\item The functions $Y_i$ {\rm [resp. $X_i$]}, $i = 0,1$, are
meromorphic in the plane $\CC_x$ cut along $[x_1 x_2] \cup [x_3 x_4]$ (resp. 
$\CC_y$ cut along $[y_1 y_2] \cup [y_3 y_4]$). In addition,
\begin{itemize}
\item $Y_0$ {\rm [resp. $X_0$]} has two zeros, no poles, and 
$\vert Y_0(x)\vert \le1, \ \forall \,|x|=1$.
\item $Y_1$ {\rm [resp. $X_1$]} has two poles and no zeros.
\item $|Y_0(x)| \leq |Y_1(x)| $ {\rm [resp. $|X_0(y)| \leq |X_1(y)|$]},
in the whole cut complex plane. Equality holds only on the cuts.  \fproof
\end{itemize}
\end{ienumerate}
\end{thm}

Combining the two basic constraints imposed on $\pi$ and $\widetilde{\pi}$ (i.e. they must be holomorphic inside their respective unit disc $\D$ and continuous on
the boundary the unit circle), and using the properties of the branches, it is possible  to make the analytic continuation of all the functions, starting from the relation
\begin{equation} \label{eq:BVP2} 
q(X_0(y), y) \pi(X_0(y)) + \widetilde{q}(X_0(y),y) \widetilde{\pi}(y)
+ \pi_0(X_0(y),y) = 0, \; y \in\CC .\end{equation}
Letting now $y$ tend successively to the upper and lower edge of the
slit $[y_1 y_2]$,  since $\widetilde{\pi}$ is
holomorphic in $\D$ and in particular on $[y_1 y_2]$, we can
eliminate $\widetilde{\pi}$ in (\ref{eq:BVP2}) to get
\[ 
\pi(X_0(y)) f(X_0(y),y) - \pi(X_1(y))f(X_1(y),y) = h(y), \; \mbox{for} \; y \in [y_1 y_2],
\]
which has exactly the profile announced in \eqref{eq:BVP}! The general theory to solve \eqref{eq:BVP} can be found in \cite{GAK,LIT2}. It involves integral forms and  an important quantity  called  the \emph{index}, defined as 
\[
\chi  \egaldef \frac{1}{2\pi} [\arg G]_\L = \frac{1}{2i\pi} [\log G]_\L,
\]
which is related to the number of existing solutions. We present now  the main substance of \cite[Theorems~5.4.1,~5.4.3]{FIM2}.

\begin{thm}\label{th:erg}
Let us introduce the following two quantities: 
\begin{equation*} 
\delta \egaldef \left\{ \begin{array}{l} 0, \; {\rm if} \; Y_0(1) < 1 \; {\rm
or} \; \left\{ Y_0(1) = 1 \; {\rm and} \; \dfrac{dq(x,Y_0(x))}
{dx}_{|x=1} > 0 \right\}, \\[0.4cm] 1, \; {\rm if} \; Y_0(1) = 1 \; {\rm
and} \; \dfrac{dq(x,Y_0(x))}{dx}_{|x=1}< 0 . \end{array}
\right.
\end{equation*}
\begin{equation*}
\widetilde{\delta} \egaldef \left\{ \begin{array}{l} 0, \; {\rm if} \;
X_0(1) < 1 \; {\rm or} \; \left\{ X_0(1) = 1 \; {\rm and} \;
\dfrac{d\widetilde{q}(X_0(y),y)}{dy}_{|y=1} > 0 \right\}, \\[0.4cm] 1,
\; {\rm if} \; X_0(1) = 1 \; {\rm and} \;
\dfrac{d\widetilde{q}(X_0(y),y)}{dy}_{|y=1}< 0. \end{array} \right.
\end{equation*}
Then  \eqref{eq:eqfonc} admits a probabilistic solution if, and only
if,
\begin{equation} \label{eq:erg}
\delta + \widetilde{\delta} = \ind{ X_0(1)=1, Y_0(1) = 1} + 1,
\end{equation}
which are the exact  conditions for the random walk to be ergodic. \fproof
\end{thm}

\begin{thm} \label{th:integral}
Under the condition \eqref{eq:erg}, the function $\pi$ is given by
\begin{equation}\label{eq:int}
\pi(x) = \frac{U(x) H(x)}{2i\pi }\int_{\M _d} {K(t)w'(t)
dt \over H^+ (t) (w(t) -w(x))} + V(x), \quad \forall x \in \mathscr{G}(\M), 
\end{equation}
where 
\begin{ienumerate}
\item $\mathscr{G}(\M)$ denotes the interior domain bounded by $\M$, and $\M_d$ is the portion of the curve $\M$ located in the lower half-plane $\Im z \leq 0$; 
\item $U,V, K$ are known functions, all involving some specific zeros of $\widetilde{q}(X_0(y),y)$ and $q(x, Y_0(x))$ inside $G_\M$; moreover $U,V$ are rational fractions;
\item $w$ is a \emph{gluing function}, which realizes the \emph{conformal mapping} of
$\mathscr{G}(\M)$ onto the complex plane cut along a segment and has
an explicit form via the Weierstrass $\wp$-function;
\item  
\begin{eqnarray*}
H(t)& =& (w(t)- X_0(y_2))^{-\widetilde{\chi}} e^{\Gamma(t)}, \; t \in
\mathscr{G}(\M), \\[0.3cm] \Gamma(t) &=& { 1 \over 2i \pi}\int_{\M _d} \log
{K(\overline{s}) \over K(s)}\, { w'(s) ds \over w(s) - w(t)} , \; t \in
\mathscr{G}(\M),\\[0.3cm] H^+(t) &=& (w(t)- X_0(y_2))^{-\widetilde{\chi}}
e^{\Gamma^+(t)},\; t \in \M_d, \\[0.3cm] \Gamma^+(t) & = & { 1 \over 2} \log
{K(\overline{t}) \over K(t)} + { 1 \over 2i \pi} \int_{\M _d} \log {
K(\overline{s}) \over K(s)}\, {w'(s) ds \over w(s) - w(t)}, \; t \in
\M_d \;.
\end{eqnarray*}
\end{ienumerate}  \fproof
\end{thm} 
The detailed proofs of these theorems can be found in the book~\cite{FIM2}.
\subsubsection{Genus $0$} $\S$ has genus $0$ if, and only if, the discriminant $D(x)$ has a multiple zero (possibly at infinity). Hence, we are left with only two branch points in the plane $\CC_x$ (resp. $\CC_y$).  This situation occurs in the $5$ following cases.

\begin{thm}\label{th6.1.1} \mbox{\ }
 The algebraic curve defined by $Q(x,y) = 0$ has genus
0 in the following cases.
\begin{enumerate}
\item $ \begin{cases} x_1  = x_2 = 0 \\ y_3  = y_4 =
\infty. \end{cases}$ \qquad
\begin{minipage}{9cm} 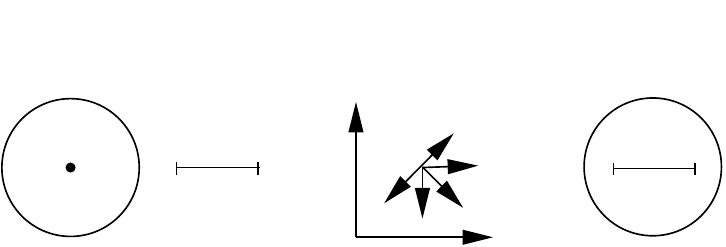 \end{minipage}
\item $ \begin{cases} y_1  = y_2 = 0 \\ x_3  = x_4 =
\infty. \end{cases}$ \qquad
\begin{minipage}{9cm} 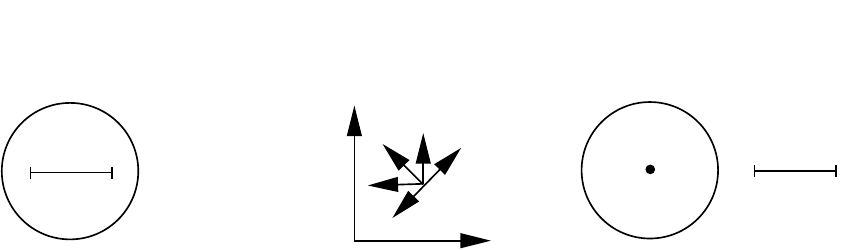 \end{minipage}
\item $ \begin{cases} x_3  = x_4 = \infty \\ y_3  = y_4 =
\infty. \end{cases}$ \quad$\;$
\begin{minipage}{9cm} 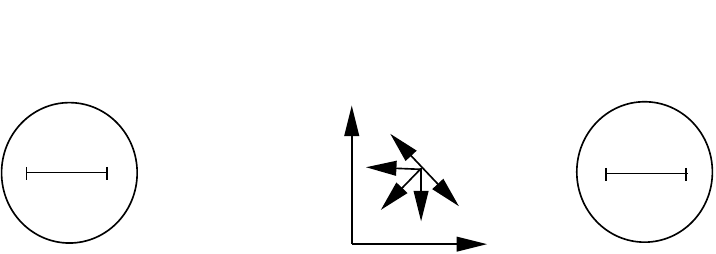 \end{minipage}
\item $ \begin{cases} x_1  = x_2 = 0 \\ y_1  = y_2 =
0. \end{cases}$ \qquad
\begin{minipage}{9cm} 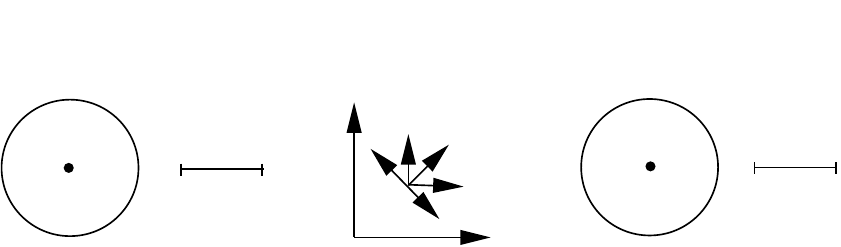 \end{minipage}
\item $ \begin{cases} x_2  = x_3 = 1 \\ y_2  = y_3 =
1. \end{cases}$ \qquad
\begin{minipage}{9cm} 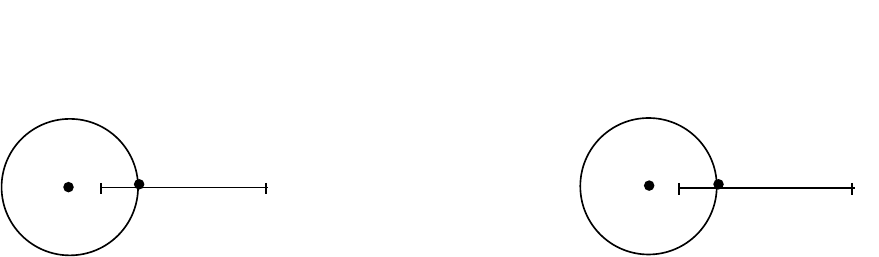 \end{minipage}
\end{enumerate}
In addition $x_3$ and $y_3$ are always positive, but $x_4$ and $y_4$
need not be positive. If for instance $x_4 < 0$, then the plane is cut
along $[- \infty, x_4] \cup [x_3, + \infty]$.
\begin{center}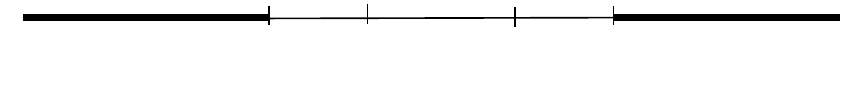 \end{center}
\begin{figure}[htb]
\begin{center} 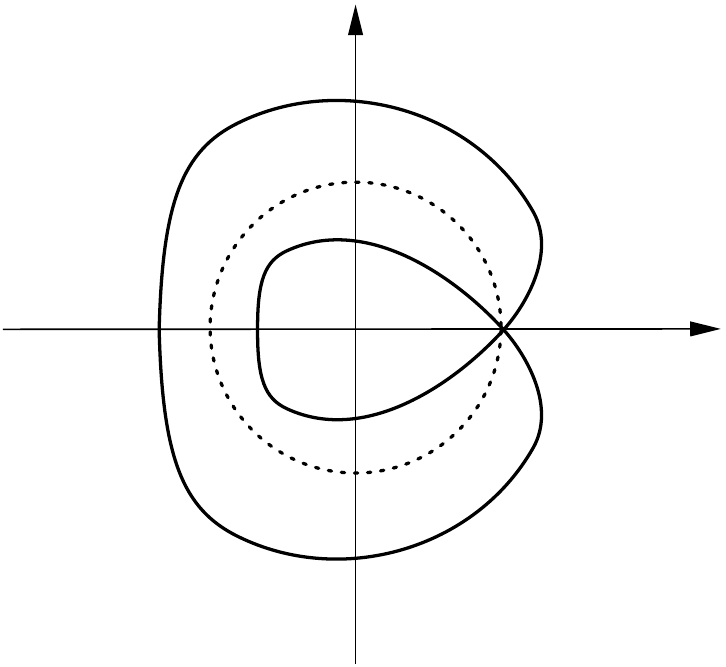 \end{center}
\vspace{-0.7cm}
\caption{The contour $\M_1\cup\M_2$, for $r<0$,}\label{fig6.6}
 $$\M_1 = X[\overrightarrow{\underleftarrow{y_1,1}}], \qquad \M_2 =
X[\overrightarrow{\underleftarrow{1,y_4}}].$$ 
\end{figure} \fproof
\end{thm}
$\bullet$ Here the algebraic curve $\{Q(x,y)=0\}$ admits of a \emph{rational uniformization} by means of rational fractions of degree $2$: that simplifies matters to a certain extent. Indeed, for each of the $5$ cases listed above,  conformal mappings (or gluing functions) can be explicitly computed, still allowing  to get integrals like in \eqref{eq:int}.  For instance, case $3$ leads to a BVP set on an \emph{ellipse}. In case $1$ (resp. case $2$) $\widetilde{\pi}(y)$ (resp. $\pi(x)$ is rational.

$\bullet$ However, case $5$ corresponds to the so-called \emph{zero drift situation} $x_2=x_3=1=y_2=y_3$ and is a bit more awkward. A  BVP can be set on the interior part  $\M_1$ of the curve shown in figure~\ref{fig6.6}, which has a corner point  at $x=1$ if, and only if, \emph{the correlation coefficient~$r$} of the random walk in the interior of the quarter plane is not zero.

\section{Examples from queueing systems} \label{sec:ex}
We shall describe the outlines of some original models using the above methods.
\subsection{Two-coupled processors (see\cite{FAIA79,FIM2})}
Consider two parallel M/M/1 queues, with infinite capacities,  under the following assumptions. \index{queueing system}
\begin{itemize}
\item Arrivals form two independent Poisson processes with parameters
$\lambda_1, \lambda_2$.
\item Service times are distributed exponentially with instantaneous service
rates $S_1$ and $S_2$ depending on the state of the system as follows.
\begin{enumerate}
\item  If both queues are busy, then $S_1=\mu_1$ and  $S_2=\mu_2$. 
\item  If queue $2$ is empty, then $S_1=\mu_1^*$.
\item  If queue $1$ is empty, then $S_2 = \mu_2^* $.
\end{enumerate}
\item The service discipline is FIFO (first-in-first-out) in each queue.
\end{itemize}
One can directly see that the evolution of the system can be described by the two-dimensional \emph{continuous time Markov process} $(M_t,N_t)$, which stands for the joint number of customers in the queues.

\smallskip
Let $p_t(m, n)$  the probability  $\Pb(M_t=m, N_t=n)$ that, at time~$t$, one finds  $m$ jobs in queue~$1$ and $n$ customers in queue~$2$. The stationary probabilities
\[
p(m,n)\egaldef \lim_{t\to\infty}p_t(m,n ) 
\]
satisfy the classical Kolmogorov forward equations, which after setting
\[
 F(x,y)= \sum_{m\ge0,n\ge0} p(m,n) x^my^n,
\]
lead to the basic functional equation (leaving the details to the reader)
\begin{equation}\label{eq:coupled-proc}
T(x,y)F(x,y) =  a(x, y)F(0, y) + b(x,y)F(x,0)  + c(x, y)F(0,0),
\end{equation}
where
\[
\begin{cases}
\DD T(x,y) = \lambda_1(1-x)+\lambda_2(1-y) 
+ \mu_1\left(1-\frac{1}{x}\right) +  \mu_2\left(1-\frac{1}{y}\right), \\[0.3cm]
\DD a(x,y)= \mu_1\left(1-\frac{1}{x}\right) + q\left(1-\frac{1}{y}\right), \\[0.3cm]
\DD b(x,y)= \mu_2\left(1-\frac{1}{y}\right) + p\left(1-\frac{1}{x}\right), \\[0.3cm]
\DD c(x,y)= p\left(\frac{1}{x}-1\right) + q\left(\frac{1}{y}-1\right), \\[0.2cm]
p= \mu_1- \mu_1^*, \\
q= \mu_2- \mu_2^*.
\end{cases}
\]
Here, there are some pleasant facts. For instance, the two roots of $T(x,y)=0$ satisfy $Y_0(x).Y_1(x)=\mu_2/ \lambda_2$, which shows that $\widetilde{\pi}(y)$ satisfies a BVP on the circle    
$\CC\bigl(\sqrt{\frac{\mu_2}{\lambda_2}}\bigr)$ with simpler formulas. 
For instance, when $pq= \mu_1\mu_2$, that is, for  $0\le\xi\le1$,
\begin{equation}\label{eq:proc-sharing}
 \begin{cases}
 \mu_1=\xi\mu_1^*,\\
 \mu_2=(1-\xi)\mu_2^*,
 \end{cases}
\end{equation}
which corresponds to the \emph{head of line processor sharing} discipline, we have, assuming the ergodicity condition $1-\lambda_1/\mu_1^*-\lambda_2/\mu_2^*>0$,
\[
F\left(0, \sqrt{\frac{\mu_2}{\lambda_2}}z\right) = 
\frac{1}{\pi} \int_0^\pi \frac{z \sin{\theta}\,v(\theta) d\theta}{z^2 -2z\cos{\theta}+1} 
+F(0,0), \quad |z|<1,
\]
where 
\begin{align*}
v(\theta) &= \frac{-\lambda_2 \sin{\theta} \,K(\theta)}
{\xi[\rho_1^*(\mu_2^*-\mu_1^*)K^2(\theta) + (\mu_1^*-\mu_2^* +\lambda_1 + \lambda_2)K(\theta) -\mu_1^*]}, \\[0.2cm]
K(\theta) & = \frac{\lambda_1+\mu_1+ \beta
 - \sqrt{[(\sqrt{(\lambda_2}+\sqrt{\mu_2})^2 +\beta][(\sqrt{\lambda_2}-\sqrt{\mu_2})^2 + \beta]}}{2\lambda_1},\\
 \beta & = \lambda_2 + \mu_2 -2\sqrt{ \lambda_2\mu_2}\cos{\theta}.
 \end{align*}
In \cite{FAIA79}, the functions $F(0,y)$ and $F(x,0)$ have been completely expressed in terms of \emph{elliptic integrals of the third kind}. 

\subsection{Sojourn time in a Jackson network with overtaking (see\cite{FIMI83,FIM2})}
\begin{figure}[htb] 
\centering 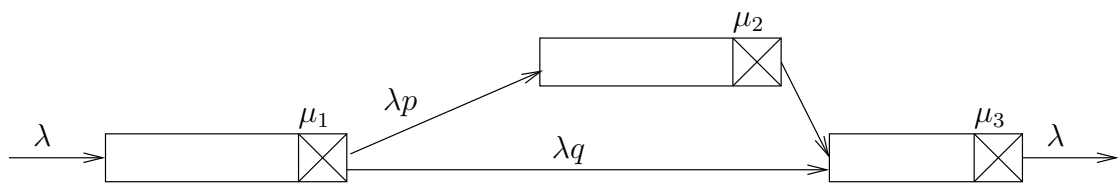
\caption{Network with overtaking} \label{fig:overtaking}
\end{figure}

A problem analyzed in \cite{FIMI83} deals with the sojourn time of a customer
in the open 3-node queueing network (of Jackson's type) shown in figure~\ref{fig:overtaking}. An inherent  \emph{overtaking phenomenon} renders things slightly more complicated. Let us just say that cutting the Gordian Knot amounts to finding the function $G(x,y,z,s)$, which is the Laplace transform of the conditional waiting time distribution of a tagged customer at a departure instant of the first queue. The following non-homogeneous functional equation can be obtained, for $|x|,|y|,|z|\le1, \Re(s)\ge0$,
\begin{align*}
K(x,y,z,s)G(x,y,z,s) &= \left(\mu_1-\dfrac{\lambda}{x}\right)G(0,y,z,s) 
+ \left(\mu_3-q\mu_1\dfrac{x}{z} - \mu_2\dfrac{y}{z}\right)G(x,y,0,s) \\[0.2cm]
&+ \dfrac{\mu_2\mu_3}{(1-x)[s+\mu_3(1-z)]},
\end{align*}
where $p,q$ are routing probabilities with $p+q=1$, $\lambda$ is the external arrival rate, 
$\mu_i$ is the service rate at queue $i$, and
\[
K(x,y,z,s)= s+ \lambda\left(1-\dfrac{1}{x}\right)+ 
\mu_1\left(1-px-q\dfrac{x}{z}\right) + \mu_2\left(1-\dfrac{y}{z}\right) +\mu_3(1-z).
\]
 Then, \emph{considering $y$ and $s$ as parameters},  the last equation takes the form
 \[
\boxed{K(x,y,z,s)\widetilde{G}(x,z)=A(x)\widetilde{G}(0,z) 
 + B(x,y,z,s)\widetilde{G}(x,0) + C(x,y,z,s)},
 \]
where $K,A,B,C$ are known functions. The reduction to a BVP is carried out according to the general methodology. Finally, setting $\rho_i=\lambda/\mu_i$, and using the geometic form  of the steady state distribution for the number of customers, it follows (see \cite{FIMI83}) that the total  sojourn time of an arbitrary customer has a Laplace transform  given by
\[
(1-\rho_1)(1-\rho_2)(1-\rho_3)\frac{\mu_1}{\mu_1+s}
G\left[\frac{\mu_1}{\mu_1+s},\rho_2,\rho_3,s\right].
\]

\subsection{Two queues with alternative service periods (see \cite{CFM87})} 
\begin{figure}[!h]
\vspace{-2cm}
\includegraphics[width=\linewidth]{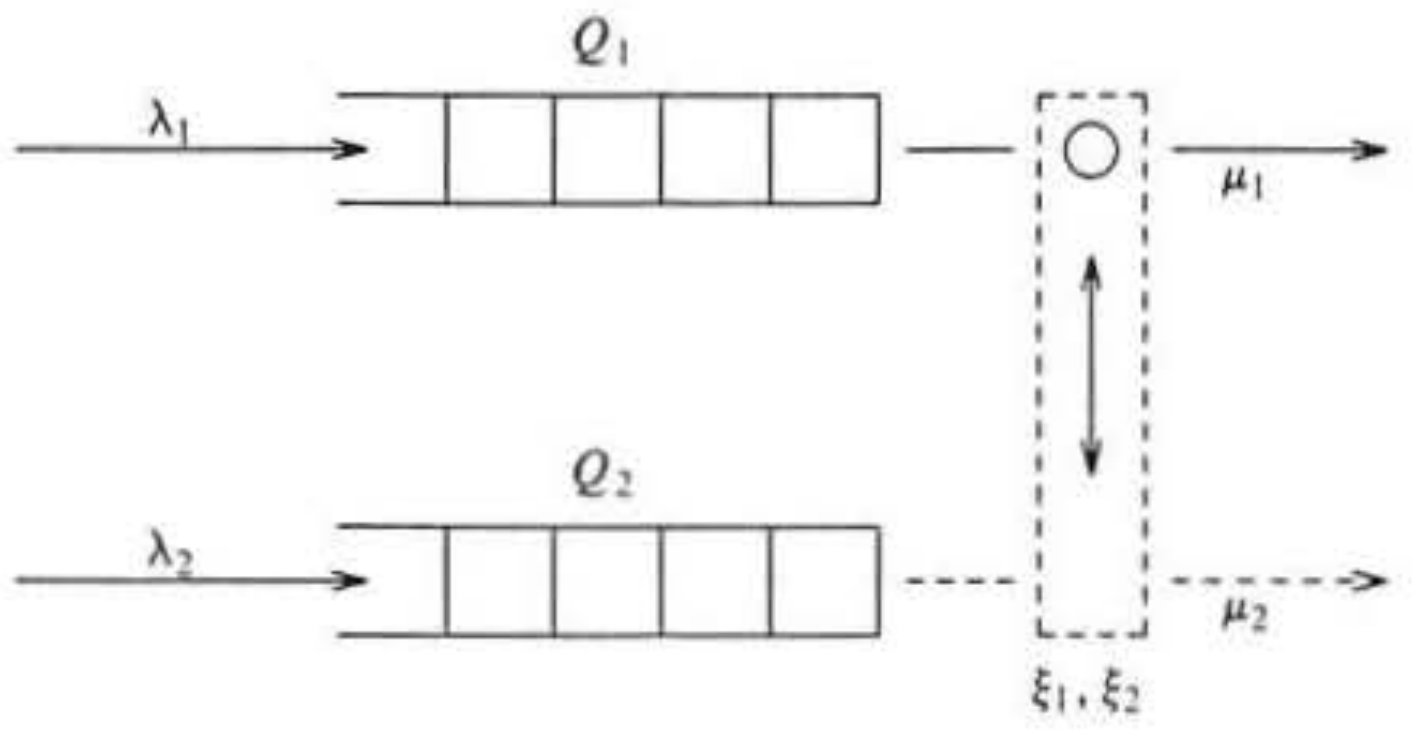}
\vspace{-14cm}
\caption{Two queues with alternative service periods}\label{fig:altserv}
\end{figure}

Consider a system of two queues, $Q_1$ and $Q_2$, and a single server that alternates service between them. (See figure~\ref{fig:altserv}) When customers are being served in $Q_i$, the system behaves as an M/M/l queue with arrival and service rate parameters $\lambda_i$ and 
$\mu_i$ ($i=1, 2$), respectively. The interarrival and service times in $Q_1$  are independent of those in $Q_2$.

Service is alternated in such a way as to limil the lime spent by the server away from a nonemply queue. If both queues are empty the server simply idles; it immediately begins service at the queue where the next arrivai occurs. When service begins at a non-empty $Q_1$, a timer is started with the initial value $T_1$. Customers in $Q_1$ are then served until eilher none remain or the $T_1$ time units have elapsed, whichever occurs first. If at this later lime, $Q_2$ is empty but $Q_1$ is still non-empty, then the above procedure is repeated; however, if $Q_2$ is non-emply then the
server begins serving customers in $Q_2$. Service of customers in $Q_2$ is similar to thal in $Q_1$; the initial timer value is now $T_2$, and a return to $Q_1$ from $Q_2$ does not occur while $Q_1$ is emply. The analysis is based on the assumption lhat $T_1$ and $T_2$  arc independcnt samples from exponential distributions wilh parameters $\xi_1$ and $\xi_2$ respcctively.

For $l=1,2$ , we define the state probabilities
\[
p_l(i,j) = \Pb(\mathrm{server \  at}\ Q_l, i\ \mathrm{customers \  in\ } Q_1, j \ 
\mathrm{customers\  in\  Q_2}),
\]
with $p(0,0) =  p_1(0, 0)=p_2(0,0)$. Then it is convenient to work with the generating functions
\[
G_1(x,y)=\sum_{i\ge1,j\ge0}p_1(i,j)x^{i-1}y^j, \quad 
G_2(x,y)=\sum_{i\ge0,j\ge1}p_2(i,j)x^iy^{j-1}.
\] 
The ergodicity  condition $\rho_1+\rho_2<1$, where $\rho_i=\lambda_i/\mu_i$, can be easily derived  by comparison with an M/G/1 queue and will be assumed to hold.

Here we end up with a \emph{system} of two FEs involving a priori \emph{four} unknown functions $G_l(0,y),G_l(x,0)$, which are easily reduced to \emph{two} by simple manipulations.

Let
\[
\begin{cases}
\DD R_1(x,y)=\lambda_1(1-x) + \lambda_2(1-y) + \mu_1\bigl(1-\frac{1}{x}\bigr), \\[0.2cm]
\DD R_2(x,y)=\lambda_1(1-x) + \lambda_2(1-y) + \mu_2\bigl(1-\frac{1}{y}\bigr), \\[0.3cm]
\DD \Delta(x,y)=(R_1(x,y)+\xi_1)(R_2(x,y)+\xi_2) - \xi_1\xi_2, \\[0.2cm]
\DD H(x)=x(R_1(x,0)+\xi_1)G_1(x,0),  \\[0.2cm]
\DD K(y) = y(R_2(0,y) + \xi_2)G_2(0,y) + \mu_2G_2(0,0) - p(0,0)\lambda_2 y, \\[0.2cm]
\DD g(x,y) = \frac{\lambda_1(1-x) + \lambda_2(1-y)}{R_2(x,y)}\xi_2 p(0,0).
\end{cases}
\]
Omitting some tiresome algebra, the following  FE of type \eqref{eq:eqfonc} can be  obtained.
\begin{lem} For $\Delta(x,y)=0$, with $|x|\le1$ and $|y|\le1$, 
we have
\begin{equation}\label{eq:alt}
H(x)-K(y)= g(x,y).
\end{equation}
\fproof
\end{lem}
Assume the polynomial $xy\Delta(x,y)$ to be irreducible. In this case, the Riemann surface $\S$ corresponding to $\Delta(x,y)=0$  is in general of \emph{genus greater than $1$}, as it reduces to a polynomial equation of degree $3$ in $x$ and in $y$ ($3$-sheeted covering). However,  there still exists a real cut, say $[x_1,x_2]$,  in the unit disc of the $\CC_x$ plane, so that a BVP of Dirichlet type (i.e. without coefficient) can be defined. Then $H(x)$ is given by a Cauchy type integral, the density of which satisfies a  Fredholm integral equation. The hassle is the analysis of the  branch points: this requires to deal with a polynomial of degree $11\,!$. Luckily enough,  a computationally  more efficient solution can be obtained via the  following approach. 
\subsubsection{Mixing uniformization and BVP}
As the uniformizalion step, we put
\begin{equation}\label{eq:unif}
R_1(x,y) + \xi_1 = \xi_1 z, \qquad R_2(x,y) + \xi_2 = \frac{\xi_2}{z}.
\end{equation}
Hence $\Delta(x,y)=0$ for $x,y$ such that \eqref{eq:unif} holds. Setting 
$\lambda=\lambda_1+\lambda_2$ and
\[
\nu=r_1x+r_2y,
\]
where $r_i=\lambda_i/\lambda$, we relate $(x, y)$ to $(\nu,z)$ by  
\begin{equation}\label{eq:x-y}
x(\nu,z) = \frac{\mu_1}{\mu_1+\lambda(1-\nu)+\xi_1(1-z)},\qquad 
y(\nu,z) = \frac{\mu_2}{\mu_2+\lambda(1-\nu)+\xi_2(1-\frac{1}{z})},
\end{equation}
with
\begin{equation}\label{eq:nu-z}
\nu = \frac{r_1\mu_1}{\mu_1+\lambda(1-\nu)+\xi_1(1-z)} + 
\frac{r_2\mu_2}{\mu_2+\lambda(1-\nu)+\xi_2(1-\frac{1}{z})}.
\end{equation}
\begin{lem}\label{lem:z-nu} \mbox{ }
\begin{ienumerate}
\item For any $z$ with $|z|=1$, there exists exactly one $\nu=\nu(z)$ such that $|v|\le1$ and $\Delta(x,y)=0$, with $|x|,|y|\le1,\,x=x(v,z),\,y= y(v,z)$. 
For $z=1$, this value is $v=1$.
\item Let $|\nu|=1$ and $\xi_1\rho_1\ge\xi_2\rho_2$. Then \eqref{eq:nu-z} has exactly one root $z(\nu)$ satisfying $|z|\ge1$; the equality $|z|=1$ holds only for $\nu=z(\nu)=1$.
\end{ienumerate} \fproof
\end{lem}
The proof of Lemma \ref{lem:z-nu} is  direct by Rouché's theorem and the principle of the argument.

\begin{figure}[!ht]
\vspace{-1.5cm}
\includegraphics[width=1.1\linewidth]{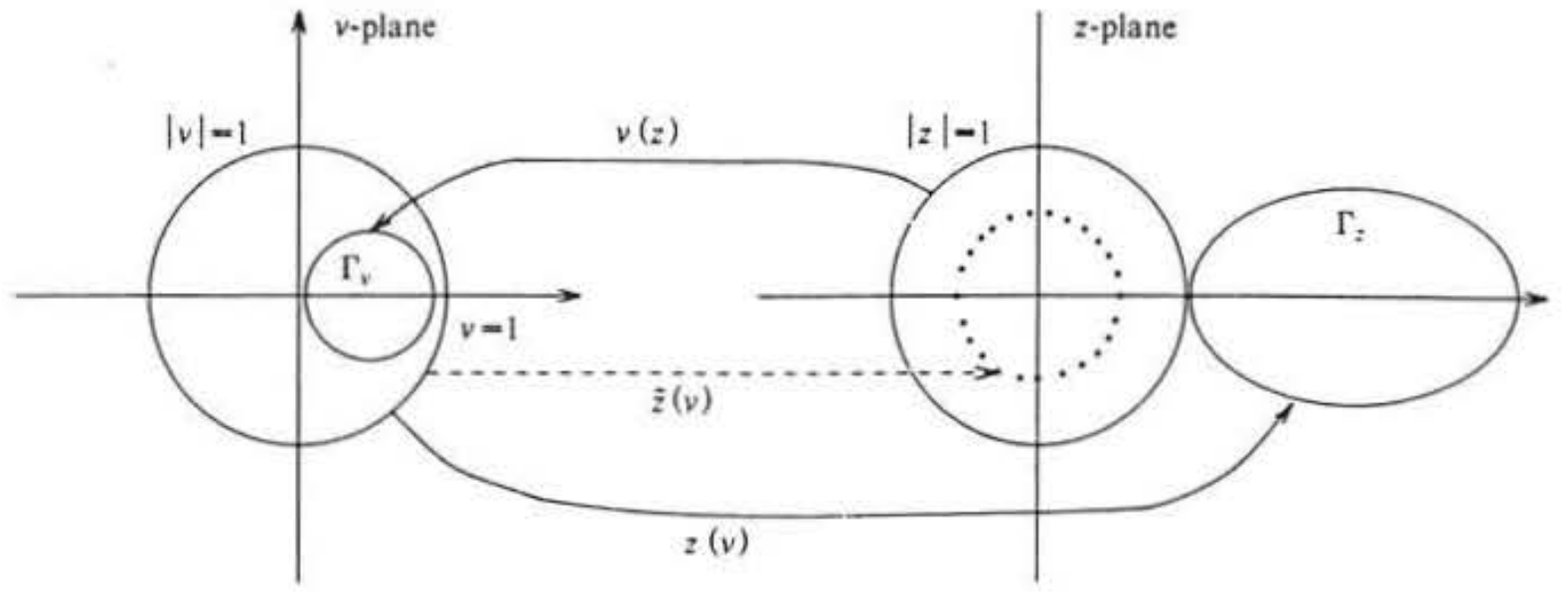}
\vspace{-17cm}
\caption{The mappings $\nu(z)$ and $z(\nu)$} \label{fig:altserv-curve1}
\end{figure}

The behaviour of the functions $\nu(z)$ and $z(\nu)$ is illustrated in 
figure~\ref{fig:altserv-curve1}. The unit circle $|z|=1$ is mapped by $\nu(z)$ onto the closed contour $\Gamma_\nu$ lying entirely within the unit circle $\nu=1$, but touching it at $\nu=1$. The unit circle $|\nu|=1$ is mapped by $z(\nu)$ onto a closed contour $\Gamma_z$ lying entirely outside and to the right of  $|z|=1$, but touching it at $z=1$. The region between  $|\nu|=1$ and  $\Gamma_\nu$ in the $\nu$-plane maps conformally onto a region outside the closed contours $|z|=1$ and $\Gamma_z$ (the latter region is not the entire region, but it does include the point al infinity). The figure also shows that the other root, say $\bar{z}(\nu)$ of \eqref{eq:nu-z} maps $|\nu|=1$ onto a contour wholly inside $|z|=1$;  this contour  does not touch $z=1$.
\begin{figure}[!ht]
\vspace{-1.5cm}
\includegraphics[width=1.1\linewidth]{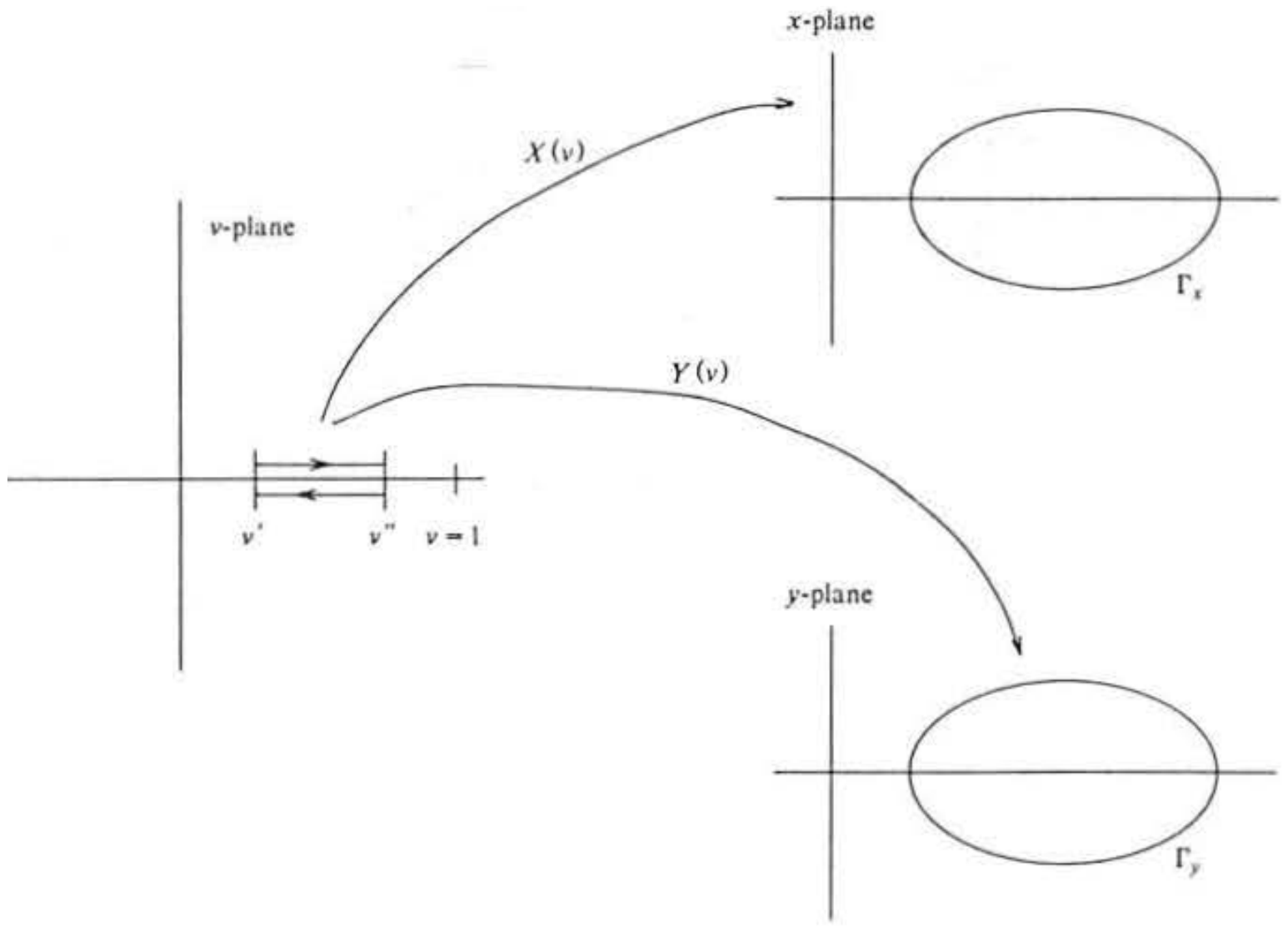}
\vspace{-13cm}
\caption{The mappings $X(\nu)$ and $Y(\nu)$} \label{fig:altserv-curve2}
\end{figure}
Writing $X(\nu)\egaldef X(\nu,z(\nu))$ and  $Y(\nu)\egaldef Y(\nu,z(\nu))$, we note that in \eqref{eq:nu-z} the desired simplification has been realized, from a cubic in $y$ (or $x$)  to a quadratic in~$z$. 

Then, $z(\nu)$ has two (real) branch points $\nu' < \nu''$ in the unit circle $|\nu|=1$. As $\nu$ moves around the cut $[\nu' < \nu'']$, $X(\nu)$ and $Y(\nu)$ traverse simple closed contours $\Gamma_x$ and $\Gamma_y$ in their respective planes $\CC_x$ and $\CC_y$, as shown in figure~\ref{fig:altserv-curve2}.  With some effort, it can also be proved that the point $y=1$ belongs to the finite region bounded by  
$\Gamma_y$, in which $K(y)$ has no pole. However, the point $x=1$  is not necessarily contained in $\Gamma_x$. Rewriting now \eqref{eq:alt} as
\begin{equation}\label{eq:alt1}
H(X(\nu))-K(Y(\nu)) = f(\nu),
\end{equation}
where $\DD f(\nu) = \frac{\lambda p(0,0)(1-\nu)z(\nu)}{1-z(\nu)}$,   $H(x)$ and $K(y)$ can be analytically continued.

The main lines of solution (which do not need  conformal mappings onto the unit disk -- see \cite{CFM87} for the details) are sketched below.

First, write the Cauchy type integrals
\[
H(x)= \frac{1}{2i\pi}\int_{\Gamma_x}\frac{H(t)dt}{t-x}, \qquad 
K(y)= \frac{1}{2i\pi}\int_{\Gamma_y}\frac{K(t)dt}{t-y}.
\]
Then, let  $x$ approach a point of the contour $\Gamma_x$, make the change of variables $t=Y(s)$, and use the Plemelj-Sokhotski formulas (see e.g., \cite{GAK}) to get
\[
\frac{1}{2}H(X(\nu)) = \frac{1}{2i\pi}\int_{L} \frac{H(X(s))X'(s)ds}{X(s)-X(\nu)}, \quad \forall \nu\in L\egaldef{[\underleftarrow{\overrightarrow{\nu',\nu''}}]}.
\]
Similarly,
\[
\frac{1}{2}K(Y(\nu)) = \frac{1}{2i\pi}\int_{\bar{L}} \frac{K(Y(s))Y'(s)ds}{Y(s)-Y(\nu)},
\quad \forall \nu\in \bar{L}\egaldef{[\underrightarrow{\overleftarrow{\nu',\nu''}}]},
\]
noting that $L$ and $\bar{L}$ are described in opposite directions. Then \eqref{eq:alt1}
leads to
\[
K(y)= \frac{1}{2i\pi}\int_{L}\frac{\psi(s)Y'(s)ds}{Y(s)-y},
\]
where $\psi(s)$ satisfies the \emph{non singular Fredholm integral equation} (which can be shown to admit a unique solution by the Fredholm alternative) 
\[
\psi(s)-\frac{1}{2i\pi}\int_L \psi(s) \frac{\partial}{\partial s} \log \frac{X(s)-X(\nu)}{Y(s)-Y(\nu)} ds = \frac{1}{2i\pi}\int_L f(s) \frac{\partial}{\partial s} \log \frac{X(s)-X(\nu)}{s-\nu}ds.
\]
Obviously,  $H(x)$ is directly obtained from \eqref{eq:alt}.

\subsection{Joining the shorter queue}\label{sec:JSQ}
This example is a long-standing problem  borrowed from queueing theory. It highlights the huge additional  complexity which arises  when space-homogeneity is only partial, even for $2$-dimensional systems. The detailed analysis can be found in \cite[Chapter~[10]{FIM2}.
\subsubsection{Equations}
Two queues with exponentially distributed service times of rate $\alpha, \beta$, respectively are placed in parallel. The external arrival process is Poisson with parameter $\lambda$. The incoming customer always joins the shorter line, or, if the lines are equal, he joins queue~$1$ or queue~$2$ with respective probabilities $\pi_1$ and $\pi_2$. The basic problem is to analyze the steady state distribution of the joint number of customers in the system. 

\smallskip
Letting  $p_t(m, n)$  denote the probability  $\Pb(M_t=m, N_t=n)$ that at time~$t$ there are  $m$ customers in queue~$1$ and $n$ customers in queue~$2$, Kolmogorov's equations for the stationary probabilities 
\[
p(m,n)\egaldef \lim_{t\to\infty}p_t(m,n ) 
\]
have to be written separately in the two distinct regions 
\[
\mathcal{R}_1\egaldef\{(m,n), m\le n\} \quad \mbox{and} \quad \mathcal{R}_2
\egaldef\{(m,n), n\le m\}.
\]
Define
\begin{equation*}\label{eq:def-JSQ}
\begin{cases}
\DD F_1(x,y)= \sum_{i,j\ge0} p(i,i+j) x^iy^j ,\quad P_1(x)= \sum_{i,\ge0} p(i,i+1) x^i, 
\\[0.5cm]
\DD F_2(x,y)= \sum_{i,j\ge0} p(i+j,i) x^iy^j ,\quad P_2(x)= \sum_{i,\ge0} p(i+1,i) x^i,
\\[0.1cm]
\DD Q(x)=F_1(x,0)=F_2(x,0)=\sum_{i,\ge0} p(i,i) x^i ,\\[0.3cm]
\DD A_1(x)=(\alpha+\lambda x)P_2(x), \quad  A_2(x)=(\beta+\lambda x)P_1 (x),\\[0.2cm]
\DD G_i(y)=F_i(0,y),\  i=1,2,\\[0.2cm]
T_1(x,y) = \lambda\left(1-\dfrac{x}{y}\right) +  \alpha\left(1-\dfrac{y}{x}\right) 
 +\beta\left(1-\dfrac{1}{y}\right),\quad R_1(x,y)= xy\,T_1(x,y),\\[0.3cm]
 T_2(x,y) = \lambda\left(1-\dfrac{x}{y}\right) +  \beta\left(1-\dfrac{y}{x}\right) 
 +\alpha\left(1-\dfrac{1}{y}\right), \quad R_2(x,y)= xy\,T_2(x,y),\\[0.3cm]
 s=\lambda + \alpha+\beta.
 \end{cases}
 \end{equation*}
Then, some direct algebra yields the following system.
\begin{equation}\label{sys:JSQ}
\begin{cases}
 T_1(x,y) F_1(x,y) & = \alpha\left(1-\dfrac{y}{x}\right)G_1(y) +
\left(\dfrac{\lambda\pi_2 y^2-\lambda x - \beta}{y}\right)Q(x) + A_1(x), \\[0.3cm]
 T_2(x,y) F_2(x,y) & = \beta\left(1-\dfrac{y}{x}\right)G_2(y) +
\left(\dfrac{\lambda\pi_1 y^2-\lambda x - \alpha}{y}\right)Q(x) + A_2(x),
\\[0.3cm]
sQ(x) & = A_1(x) + A_2(x).
\end{cases}
\end{equation}
\subsubsection{Reduction of the number of unknown functions} 
The ergodicity of the process  is equivalent to the existence of $F_1(x,y)$ and 
$F_2(x,y)$ holomorphic in $\D\times\D$ and continuous in 
\mbox{$\overline{\D}\times\overline{\D}$}. At first sight, system \eqref{sys:JSQ} includes four unknown functions of one variable. In fact, this number immediately  boils down to two,  by using the zeros of the \emph{kernels} \mbox{$R_j(x,y),j=1,2,$} in  
\mbox{$\overline{\D}\times\overline{\D}$}. So,  we are  left with \emph{two unknown functions of one complex variable}, say for instance \mbox{$G_i(\cdot), i=1,2$}, or \mbox{$A_i(.), i=1,2$}. In order to get additional information, we can  combine  several BVPs  derived from system \eqref{sys:JSQ}. As we shall see, determining one function is sufficient to find all the others.
\subsubsection{Meromorphic continuation to the complex plane}\label{sec:MeroCont}
It is easy to check that the algebraic curves defined by \mbox{$\{R_j(x,y)=0\}, j=1,2$}, correspond to random walks of genus~$0$, case~$3$  of  Theorem~\ref{th6.1.1}.

\paragraph{Notation and Assumption} \emph{For convenience and to distinguish between the two kernels, we shall add, either in a superscript or subscript position ad libitum, the pair $\alpha\beta$ (resp. $\beta\alpha$) to any quantity related to the kernel $R_1(x,y)$ (resp. $R_2(x,y)$). For instance, the branches  $Y_0^\ab, X_1^\ab$, etc. Also, if a property holds both for $\ab$ and $\ba$, the pair is omitted.  \textbf{From now on, we assume \mbox{$\beta>\alpha$}}, the case  \mbox{$\beta=\alpha$} being considered in a separate section}.

The functions $Y_i^\ab(x), i=0,1$, have exactly two branch points, which are always located  inside $\D$,
\begin{equation}\label{eq:bpY}
x=0, \quad \mathrm{and} \quad x_\ab^* = \frac{4\alpha\beta}{s^2-4\alpha\lambda}<1.
\end{equation}
With the notation of Section~\ref{sec:genus1} for the contour corresponding to a slit, 
$\Psi^\ab$ will denote the ellipse obtained  by the mapping
\[
[\underrightarrow{\overleftarrow{0,x^*_\ab}}] \ \overset{Y^\ab}{\longrightarrow} \ \Psi^\ab. 
\] 
Note that $\beta>\alpha \iff x_\ab^*< x_\ba^*$.
In the $\CC_y$-plane, setting $y=u+iv$, the equation of $ \Psi^\ab$ is
\[
\left(u- \frac{\beta s}{s^2-4\alpha\lambda}\right)^2 + \frac{s^2v^2}{s^2-4\alpha\lambda} = 
\frac{\beta^2s^2}{(s^2-4\alpha\lambda)^2}.
\]
Similarly for $X_i^\ab(x), i=0,1$, with the branch-points
\[
y_1^\ab = \dfrac{\beta}{s+2\sqrt{\alpha\lambda}} \ < \ 
y_2^\ab = \dfrac{\beta}{s-2\sqrt{\alpha\lambda}}, \quad   0<y_1^\ab<y_2^\ab<1,
\]
and $\Phi^\ab$ will denote the ellipse obtained by the mapping
\[
[\underrightarrow{\overleftarrow{y_1^\ab,y_2^\ab}}]
 \ \overset{X^\ab}{\longrightarrow} \ \Phi^\ab.
\]
In the $\CC_x$-plane, setting $x=u+iv$,  the equation of $\Phi^\ab$ is
\[
\left(u- \frac{2\alpha\beta}{s^2-4\alpha\lambda}\right)^2 + \frac{s^2v^2}{s^2-4\alpha\lambda} = \frac{\alpha\beta^2s^2}{\lambda(s^2-4\alpha\lambda)^2}.
\]
Exchanging  the parameters $\alpha$ and $\beta$, the respective branch points of  
$Y_i^\ba(x)$ are
\begin{equation}\label{eq:bpX}
x=0, \quad \mathrm{and} \quad x_\ba^* = \frac{4\alpha\beta}{s^2-4\beta\lambda}<1,
\end{equation}
and those of $X_i^\ba(x)$,
\[
y_1^\ba = \dfrac{\alpha}{s+2\sqrt{\beta\lambda}} \  <
y_2^\ba = \dfrac{\alpha}{s-2\sqrt{\beta\lambda}}, \quad 0<y_1^\ba<y_2^\ba<1.
\]

\begin{thm}\label{theo:AnalCont}
The functions $Q(\cdot),G_i(\cdot),A_i(\cdot), i=1,2,$ can be continued as meromorphic functions to the whole complex plane. \fproof
\end{thm}
The proof  (first established in 1979) can be found in \cite{FIM2} and relies on the following lemma.
\begin{lem} \label{lem:Dn} Let $\D_n$ be the domain recursively defined by
\[
\begin{cases}
\D_0  & \!\!= \ \D, \\
\D_{n+1} & \!\!= \inf{\left\{(X_1\circ Y_1)^\ab(\D_n), (X_1\circ Y_1)^\ba(\D_n)\right\}}.
\end{cases}
\]
Then $\D_n \subset \D_{n+1}$ and $\DD \lim_{n\to\infty}\D_n = \CC$, where 
$\CC$ denotes the complex plane. \fproof
\end{lem}
\subsubsection{Functional equation for $G_1(y)$ and integral equation for $Q(x)$}
Hereafter, we shall list the main results of this study, presented in the form of a global Proposition. Proofs involve sharp technicalities are omitted. They can be found in~\cite{FIM2} and references therein.

Let
\[
\begin{cases}
\Delta(x) &= s- \left(\frac{\alpha}{x} + \lambda\pi_2\right)Y_0^\ab(x) -  \left(\frac{\beta}{x} + \lambda\pi_1\right)Y_0^\ba(x), \\[0.3cm]
K_1(z) &= \dfrac{\biggl(1-\dfrac{Y_0^\ab(z)}{z}\biggr)\Delta_1(z)}
{\biggl(1-\dfrac{Y_0^\ba(z)}{z}\biggr) W^\ab(z)}, \quad 
L_1(z) = \dfrac{\biggl(1-\dfrac{Y_1^\ab(z)}{z}\biggr)\Delta(z)}
{\biggl(1-\dfrac{Y_0^\ba(z)}{z}\biggr) W^\ab(z)}, \\[1cm]
\Delta_1(z) &= s - \left(\dfrac{\alpha}{z} + \lambda\pi_2\right)Y_1^\ab(z) -  
\left(\dfrac{\beta}{z} + \lambda\pi_1\right)Y_0^\ba(z).
\end{cases}
\]

\begin{prop}\label{theo:G1-Q} The two functions $G_1(y)$ and $Q(x)$ have the following properties. 
\begin{enumerate}
\item 
\begin{equation}
\boxed{\begin{aligned}\label{eq:BVP-G1} 
\lefteqn{L_1(z)G_1(Y_1^\ab(z)) - L_1(\bar{z})G_1(Y_1^\ab(\bar{z})) \ =} \\
 & K_1(z)G_1(Y_0^\ab(z)) - K_1(\bar{z})G_1(Y_0^\ba(\bar{z})), \ z\in\Phi^\ba
\end{aligned}}\,,
\end{equation}
which is equivalent to a  generalized Riemann-Carleman BVP on the closed contour 
$\L= Y_1^\ab(\Phi^\ba)$, having a unique solution analytic in $\D$ if and only 
if 
\[
\Delta'(1)>0, \quad \mbox{\emph{or equivalently}} \  \lambda<\alpha+\beta.
\]
Moreover, under the above ergodicity condition, the whole system~\eqref{sys:JSQ} has also an analytic solution in $\D$. 
\item Equation \eqref{eq:BVP-G1} reduces to a BVP of the form
\begin{equation} \label{eq:BVP-U}
\boxed{ U^+(t) - U^-(t) = H(t)\overline{U^+(t)} + C, \quad t\in w(\Phi^\ab)}\,,
\end{equation}
where $H(t)$ is known, and $w(z)$ is a conformal gluing of the domain inside the ellipse $\Phi^\ab$ onto the complex plane cut along an open smooth arc. Moreover, \eqref{eq:BVP-U} defines a \emph{Noetherian operator with index $0$}, that is a Fredholm operator.
\item For $y\in\CC_y$, the function $G_1(y)$  satisfies the non local FE
\begin{equation}
\boxed{\begin{aligned}\label{eq:Merom-G1}
\lefteqn{L_1(X_0^\ba(y))G_1(Y_1^\ab\circ X_0^\ba(y)) - L_1(X_1^\ba(y))
G_1(Y_1^\ab\circ X_1^\ba(y))\ =}\\
 & K_1(X_0^\ba(y))G_1(Y_0^\ab\circ X_0^\ba(y)) -
  K_1(X_1^\ba(y))G_1(Y_0^\ab\circ X_1^\ba(y))
\end{aligned}}\,.
\end{equation}
\item The function $Q(x)$ satisfies the real  integral equation
\begin{equation} \label{eq:Q}
\boxed{Z(x)Q(x) = \int_0^{x^*_\ba}Q(u)S(x,u)du + W(x), \quad \forall x\in [0,x^*_\ba]}\,,
\end{equation}
where $Z(x), S(x,u), W(x)$ are known quantities. \fproof
\end{enumerate}
\end{prop}
Some facts have to be stressed.
\begin{itemize}
\item All the unknown functions appearing in system~\eqref{sys:JSQ} can be determined once  $Q(.)$ is known.
\item From \eqref{eq:Merom-G1}, explicit (but intricate) recursive computations of the poles and residues of $G_1(\cdot)$ can be achieved (see  \cite{COHE98}).
\item The  integral equation \eqref{eq:Q} \emph{contains both a regular part and a singular part}.
\end{itemize}
\subsection{Explicit integral forms for equal service rates ($\alpha=\beta$)}
In this case, the problem simplifies in a breathtaking way\,! Indeed, the two kernels $T_i(x,y),i=1,2,$ are equal. Using the same objects as before (and omitting the indices 
$\ab$ or $\ba$), we put
\[
\begin{cases}
F(x,y)  = F_1(x,y) + F_2(x,y), \\[0.2cm]
G(x,y) = G_1(x,y) + G_2(x,y), \\[0.2cm]
\Delta(x,y)  = s-\left(\dfrac{2\alpha}{x}+\lambda\right)y,
\end{cases}
\]
Then, from system \eqref{sys:JSQ}, we get the reduced functional equation
\begin{equation} \label{sys:JSQ1}
F(x,y) T(x,y) =  \alpha\left(1-\dfrac{y}{x}\right)G(y) -  \Delta(x,y)Q(x).
\end{equation}
 the resolution of which \eqref{sys:JSQ1} is straightforward by applying the methods previously discussed. The  result is presented in the following proposition 
 without further comment. 
\begin{prop} \mbox{ }
When $\alpha=\beta$, the system is ergodic if and only if  $\lambda<2\alpha$, and in this case
\begin{equation}\label{eq:a=b}
Q(x) = K\,e^{\Gamma(x)},
\end{equation}
\end{prop}
\[
\begin{cases}
\Gamma(x) &= \DD \dfrac{1}{2i\pi} \int_\Phi \dfrac{\log{w(t)}\theta'(t)dt}
{\theta(t)- \theta(x)}, \quad x\in \mathscr{G}(\Phi) ,\\[0.5cm]
w(t) &= \dfrac{\overline{i\xi(t)}}{i\xi(t)},\quad \xi(t)= \dfrac{\Delta(t)}{1-\frac{Y_0(t)}{t}},
\quad \Delta(t) =   \Delta(t,Y_0(t)),
\end{cases}
\] 
where $K$ is a positive constant, $\mathscr{G}(\Phi)$ is the interior domain bounded by the ellipse $\Phi$, and $\theta(.)$ denotes the conformal mapping of 
$\mathscr{G}(\Phi)$ onto the unit disc.

\section{Counting lattice walks in the quarter plane}\label{sec:counting}
Enumeration of planar lattice walks has become a classical topic in combinatorics. For 
a given set $\Sc$ of allowed jumps (or steps), it is a matter of counting the number of paths starting from some point and ending at some arbitrary point in a given time, and possibly restricted to some regions of the plane.

Then three important questions naturally arise.
\begin{itemize}
     \item[Q1:] \label{q-howmany} How many such paths exist?
     \item[Q2:] \label{q-nature} What is the nature of the associated \emph{counting generating function} (CGF)  of the numbers of walks? Is it \emph{holonomic},\footnote{A function of several complex variables is said to be holonomic if the vector space over the field of rational functions spanned by the set of all derivatives is finite dimensional. In the case of one variable, this is tantamount to saying that the function satisfies a linear differential equation where the coefficients are rational functions (see e.g., \cite{FlSe}).} and, in that case, \emph{algebraic} or even  \emph{rational}?  \index{holonomic}
  \item[Q3:] \label{q-asymptotic} What is the asymptotic behavior, as their length goes to infinity, of the number of walks ending at some given point or domain (for instance one axis)? 
\end{itemize}
If the paths are not restricted to a region, or if they are constrained to remain in a half-plane,  the CGFs  have an explicit form and can only be rational or algebraic (see \cite{BMP2003}). The situation happens to be much richer if the walks are confined to the quarter plane $\Zb_+^2$.

So, we shall focus on walks confined to $\mathbb{Z}_+^2$, starting at the origin and having \emph{small steps}. This means exactly that  the set $\mathcal S$ of admissible steps is included in the  set of the eight nearest neighbors, i.e., $\mathcal S \subset \{-1, 0, 1\}^2 \setminus\{(0, 0)\}$. By using  an extended Kronecker's delta, we shall write
\begin{equation}
\label{def_deltaij}
     \delta_{i,j}=\left\{\begin{array}{ccc}
     1 &\text{if} & (i,j)\in\mathcal S,\\
     0 &\text{if} & (i,j)\notin\mathcal S.
     \end{array}\right.
\end{equation}
A priori, there are $2^{8}$ such models. In fact, after eliminating trivial  cases and models equivalent to walks confined to a half-plane, and noting also that some models are obtained from others by symmetry, it was shown in \cite{BMM2010} that one is left with  $79$ inherently different problems to analyze.

A common starting point to deal with these $79$ walks relies on the following analytic approach. Let $f(i,j,k)$ denote the number of paths in $\Zb_{+}^{2}$ starting from $(0,0)$ 
and ending at $(i,j)$ at time $k$ (or after $k$ steps). Then the corresponding CGF 
     \begin{equation}\label{eq:tri}
          F(x,y,z)=\sum_{i,j,k\geq 0}f(i,j,k)x^{i}y^{j}z^{k}
     \end{equation}
satisfies the functional equation (see \cite{BMM2010} for the details)
     \begin{equation}  \label{eq:cgf}
         \boxed{ K(x,y,z)F(x,y,z) =  c(x)F(x,0,z)+\widetilde{c}(y)F(0,y,z) + c_0(x,y,z)},
     \end{equation}
valid a priori in the domain $|x|\leq1,|y|\leq1,|z|<1/|\mathcal{S}|$, where 
\begin{equation*}
\begin{cases}
\DD K(x,y;z)  =  xy\Bigg[\sum_{(i,j)\in\Sc} x^iy^j-1/z\Bigg] , \\[4ex]
\DD c(x)  = \sum_{i \leq 1} \delta_{i,-1}  x^{i+1}  ,\quad
\DD \widetilde{c}(y)=\sum_{j \leq 1} \delta_{-1,j} y^{j+1}, \\[0.cm]
\DD c_0(x,y,z) = -\delta_{-1,-1} F(0,0,z)-x y/z ,
\end{cases} 
\end{equation*}
\subsection{Goup classification of the 79 main random walks}
For various reasons (mathematical, but also related to computational efficiency), it seems of interest to get information about the nature of the generating functions. In the probabilistic context  of  Sections~\ref{sec:2var}-\ref{sec:ex}, assuming the group $\H$  to be finite (see definition \eqref{eq:gr2n}), first results in that direction have been given
 in~\cite{FIM2} in terms of necessary and sufficient conditons  for the unknown functions to be algebraic or rational. 

We  write $h_{\alpha} \egaldef
\alpha (h)$, for all automorphisms $\alpha$ and all functions $h$
belonging to $\CC_Q(x,y)$. For any $h \in\CC_Q(x,y)$, let the \emph{norm} $N(h)$ be defined as
\begin{equation}\label{eq-norm}
N(h)  \egaldef  \prod_{i=0}^{n-1} h_{\delta^i}. 
 \end{equation}

Written on $Q(x,y)=0$, equation \eqref{eq:eqfonc} yields the system
\begin{equation} \label{eq.pi}
 \begin{cases}
 \pi= \pi_\xi, \\[0.1cm]
 \pi_\delta - f \pi = \psi,
\end{cases}
 \end{equation}
where 
\begin{equation*}
f= \frac{q\widetilde{q}_\eta}{\widetilde{q}q_\eta} , \quad
 \psi = \frac{q_0\widetilde{q}_\eta}{q_\eta \widetilde{q}} - \frac{(q_0)_\eta}{q_\eta}.
\end{equation*}
It has been shown in \cite{FR2010} (see Theorem~11.3.3 in~\cite{FIM2}) that, when $n$ is finite and $N(f)=1$, the solution $\pi(x,y)$ of \eqref{eq:eqfonc} is always \emph{holonomic}.

In \cite{BMM2010}, the authors  consider the group 
$W\egaldef\langle \alpha,\beta\rangle$ generated by the two  birational transformations leaving invariant the generating function $\sum_{(i,j)\in\Sc}x^{i}y^{j}$,
     \begin{equation*} 
          \alpha(x,y)= \Bigg(x,\frac{1}{y}\frac{\sum_{(i,-1)\in\Sc}x^{i}}
          {\sum_{(i,+1)\in\mathcal{S}}x^{i}}\Bigg),\qquad
          \beta(x,y)=\Bigg(\frac{1}{x}\frac{\sum_{(-1,j)\in\mathcal{S}}y^{j}}
          {\sum_{(+1,j)\in\mathcal{S}}y^{j}},y\Bigg).
     \end{equation*}

Clearly $\alpha^2=\beta^2=\text{Id}$, and $W$ is a dihedral group of even order (always~$\ge4$).
 
\begin{figure}[!ht]
\hfill\hbox{}
\includegraphics[width=2.5cm]{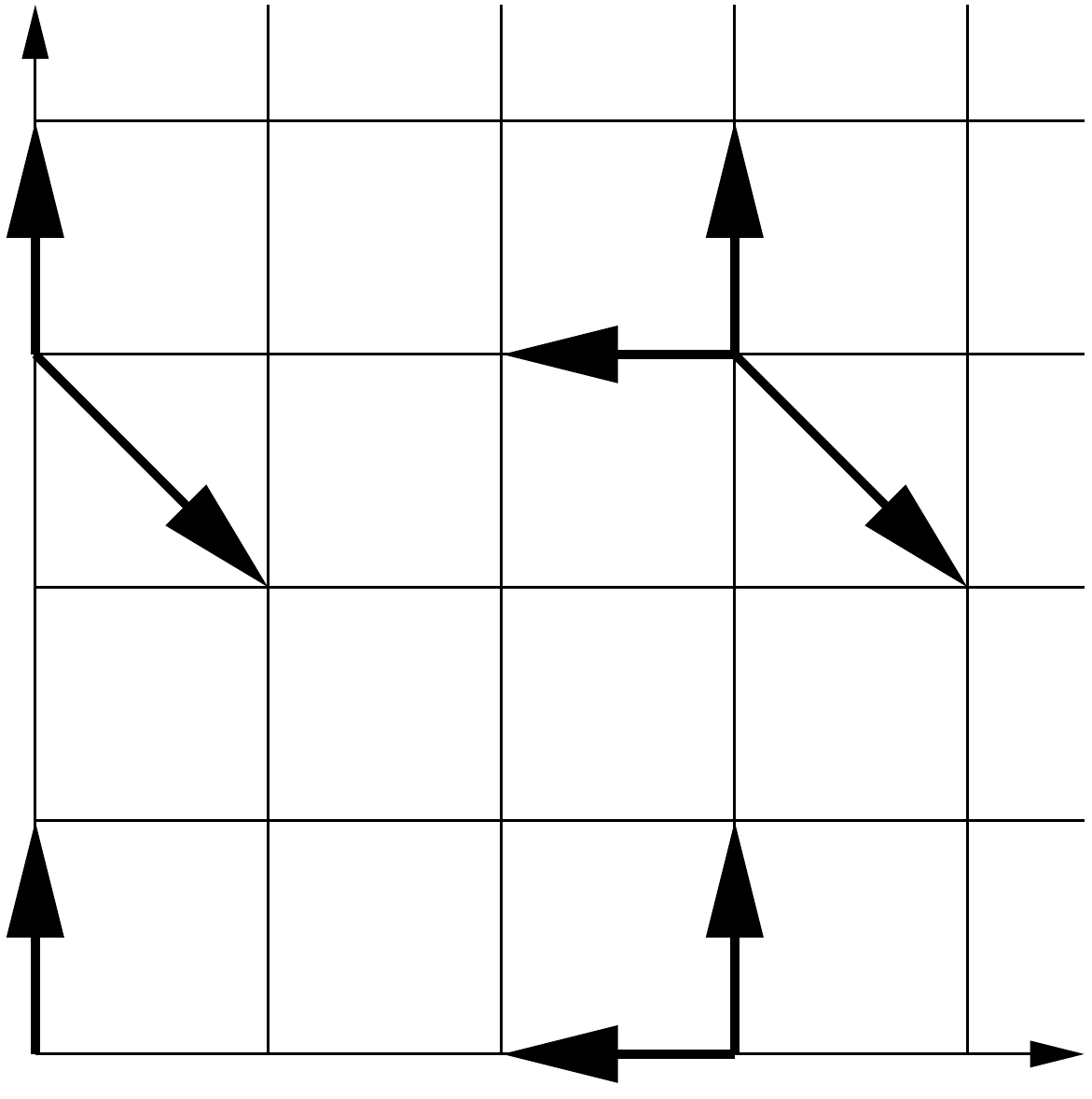}
\hfill\hbox{}
\includegraphics[width=2.5cm]{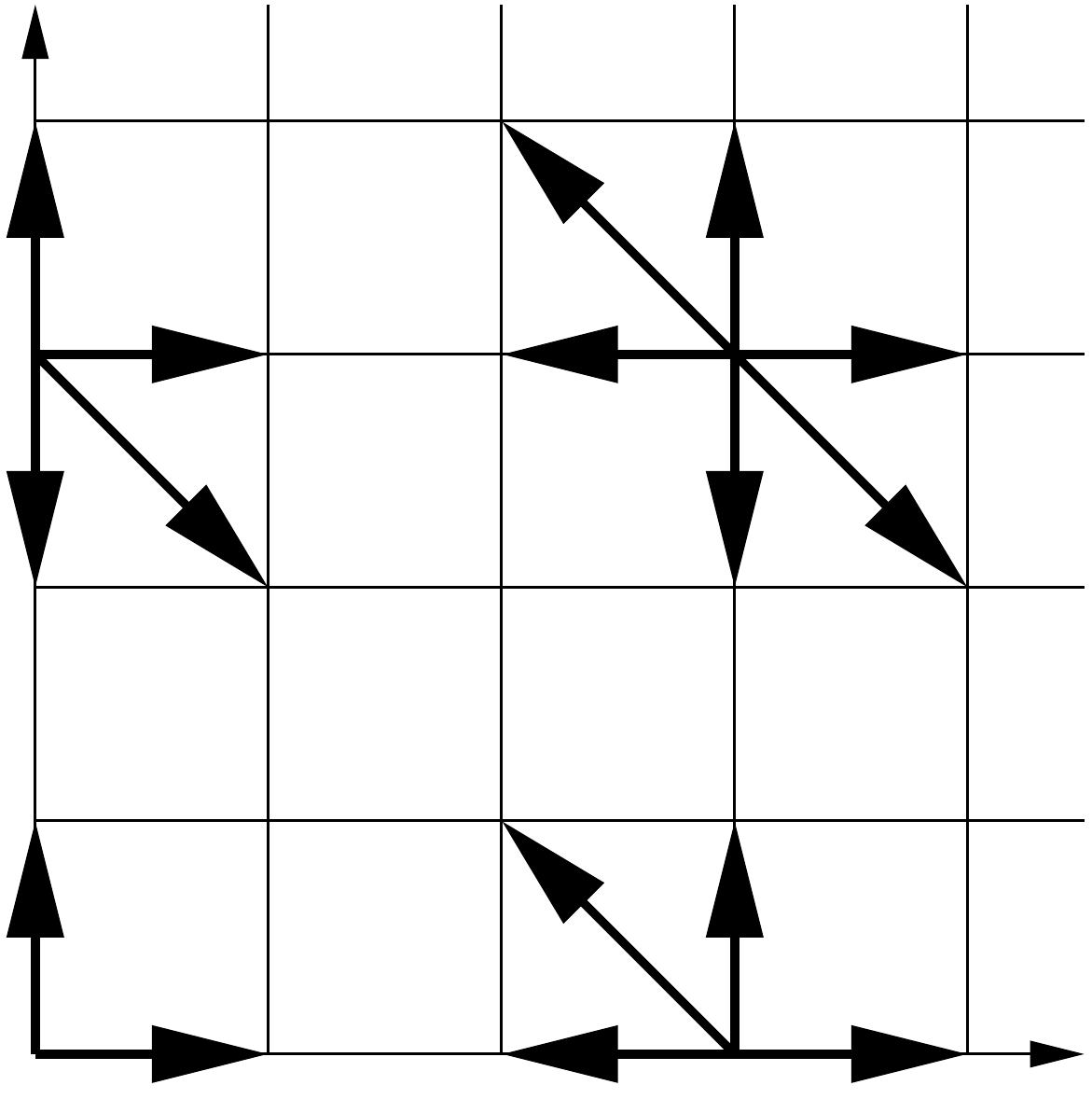}
\hfill\hbox{}
\includegraphics[width=2.5cm]{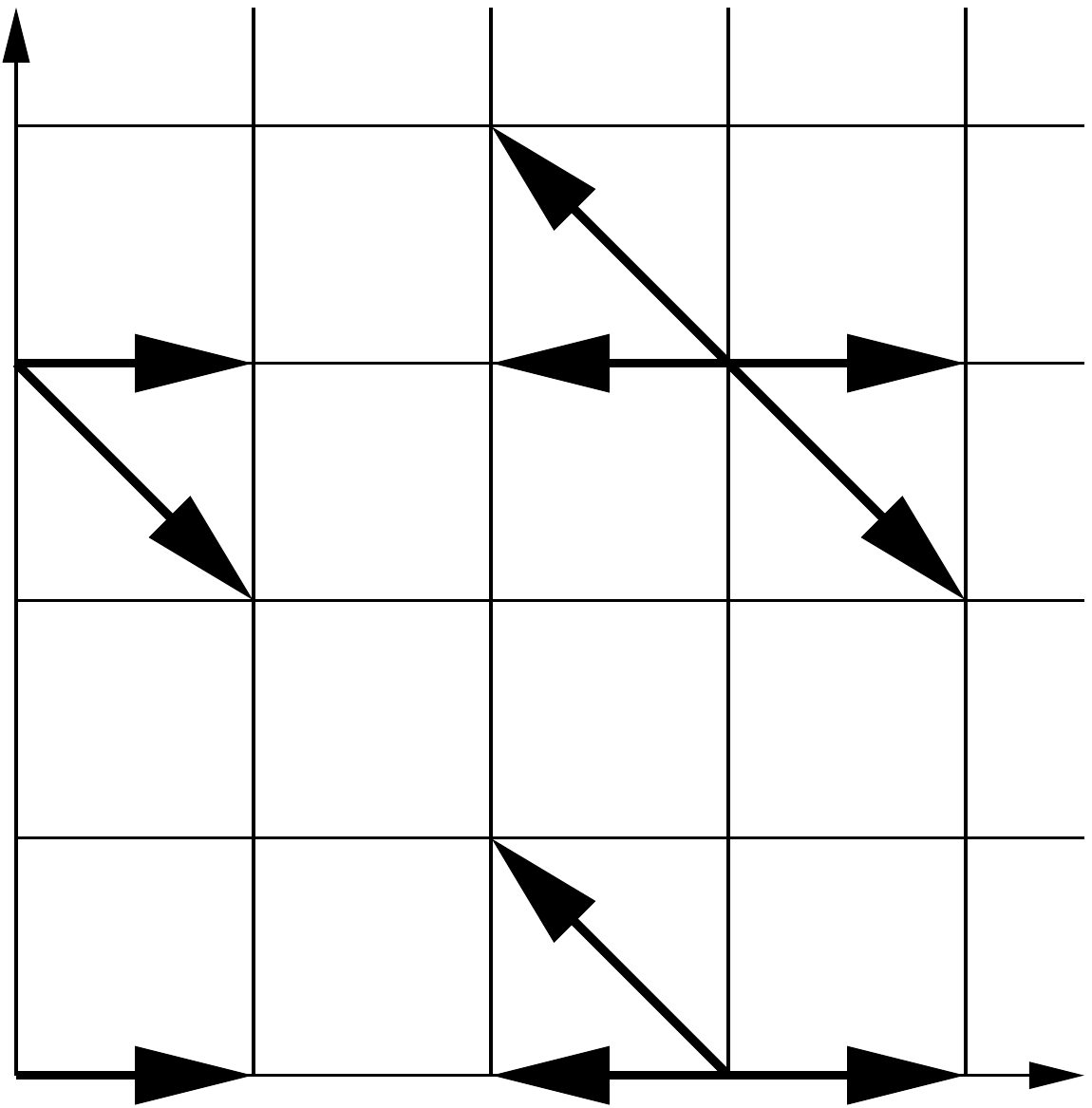}
\hfill\hbox{}
\caption{On the left 2 walks with a group of order 6. On the right, 1 walk with a group of order 8.} \label{Fig_hol}
\end{figure}

\begin{figure}[!ht]
\hfill\hbox{}
\includegraphics[width=2.5cm]{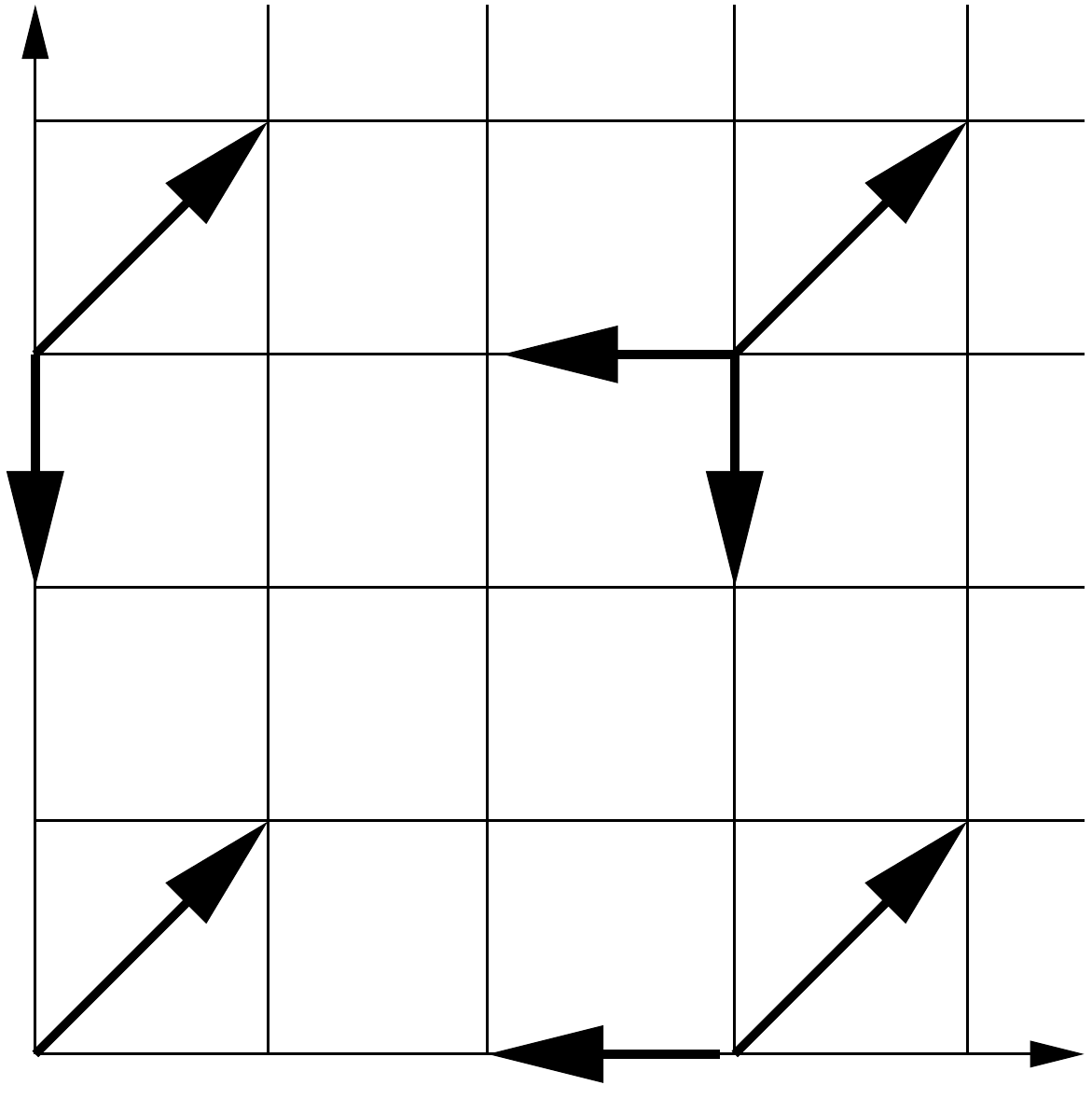}
\hfill\hbox{}
\includegraphics[width=2.5cm]{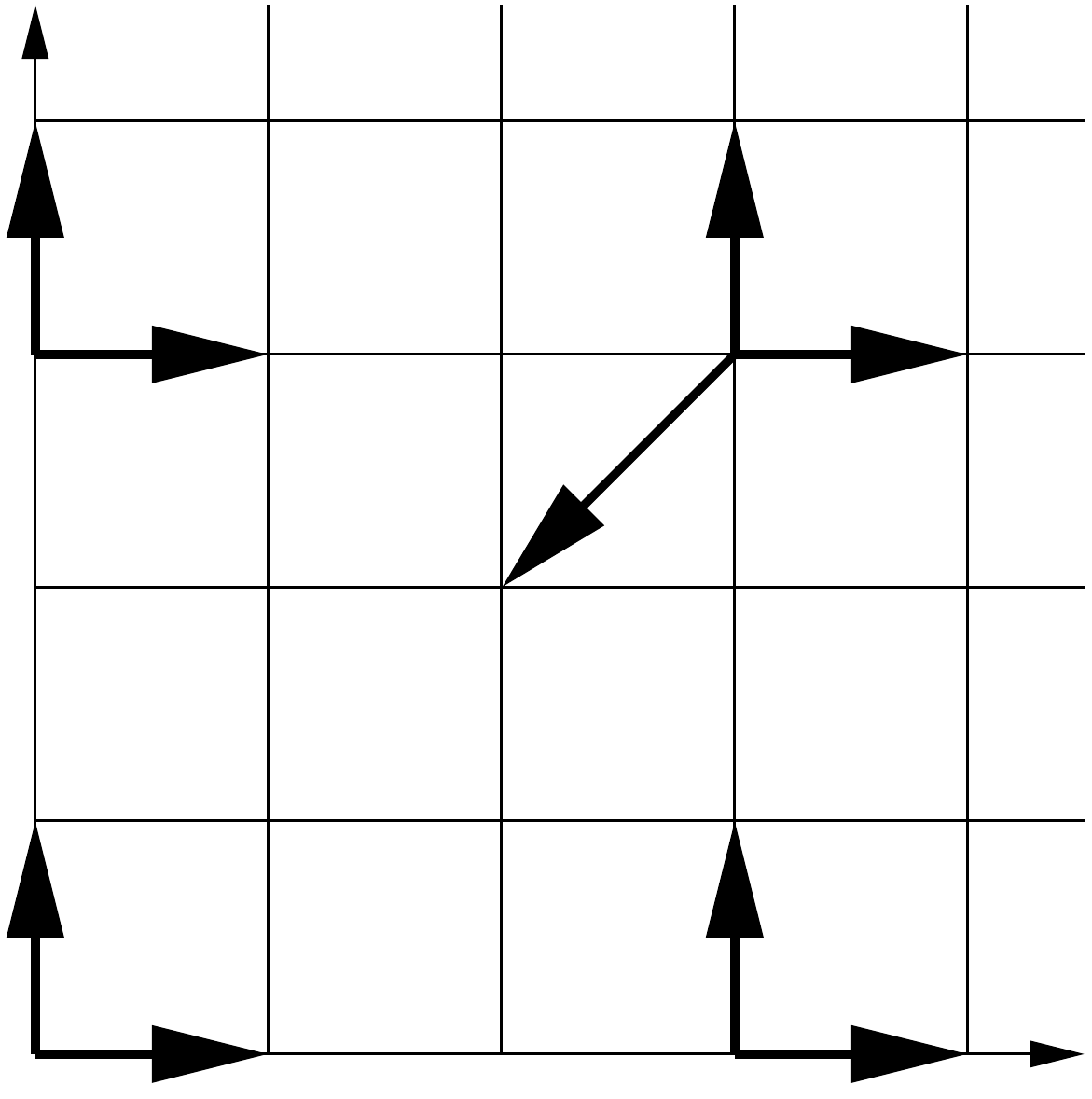}
\hfill\hbox{}
\includegraphics[width=2.5cm]{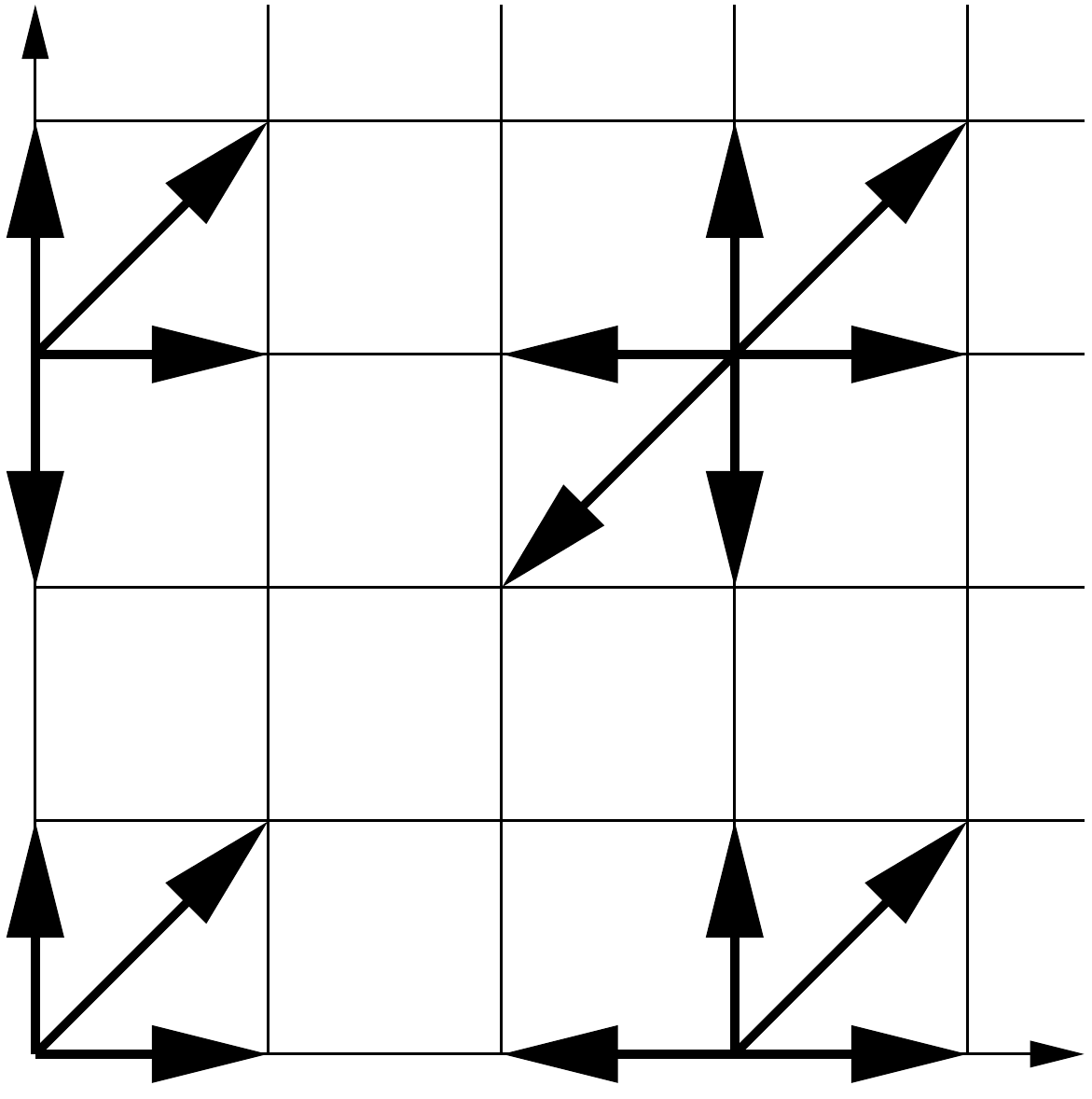}
\hfill\hbox{}
\includegraphics[width=2.5cm]{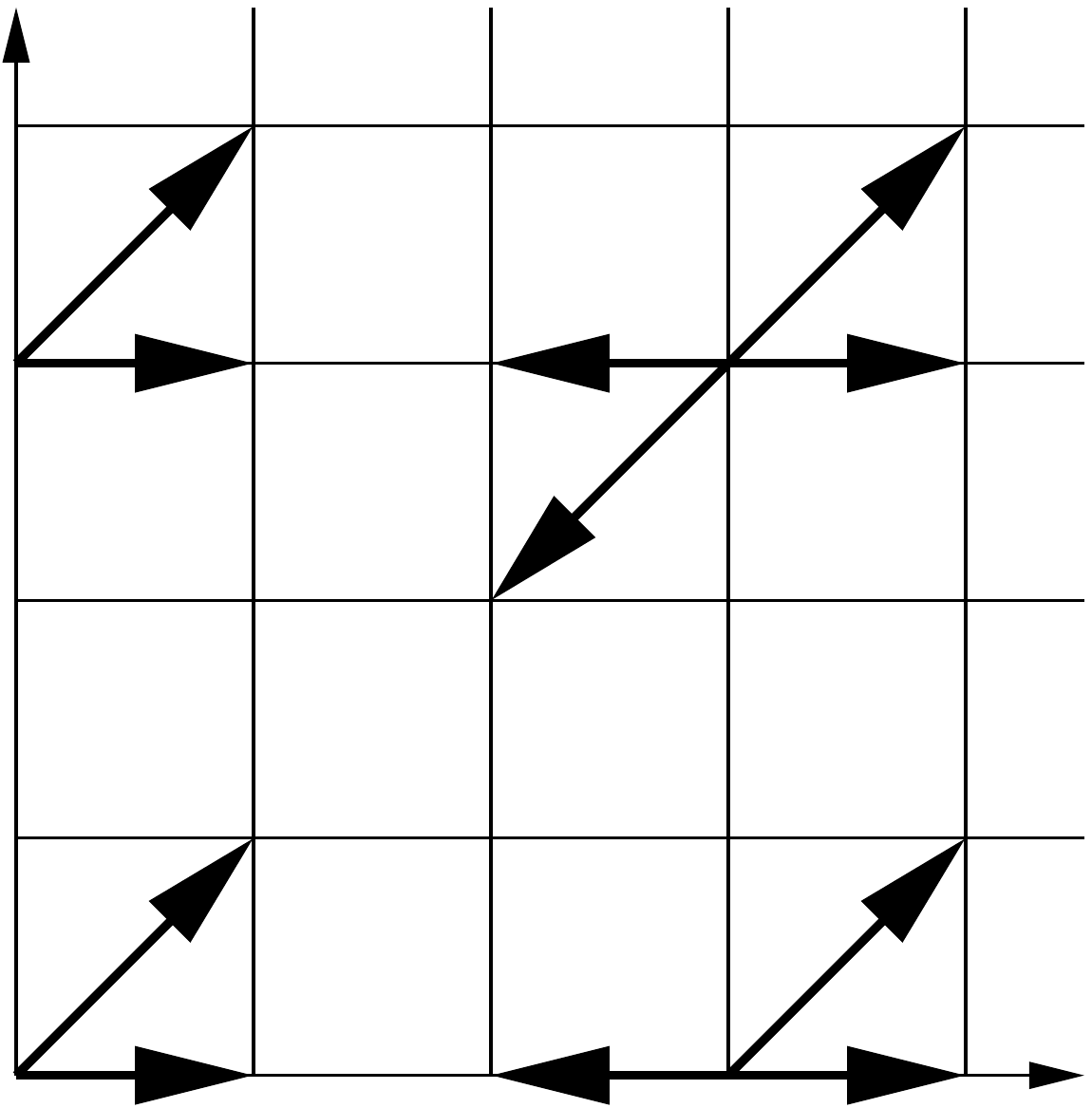}
\hfill\hbox{}
\caption{On the left,  3 walks with a group of order 6. 
         On the right, 1 walk with a group of order 8.} \label{Fig_alg}
\end{figure}

 The difference between the groups $W$ and $\H$ (see  Definition~\ref{def:gr}) is not only of a formal character. In fact, $W$ is defined on all of $\CC^2$, whereas $\H$ acts only on the algebraic curve defined by of the type (see \eqref{eq:FE}). Clearly 
\begin{equation}\label{HW}
\textit{Order}(\H)\le \textit{Order}(W),
\end{equation}
and a quick analysis shows that the group $W$ is, in some sense, less general than 
$\H$, since it must keep $\sum_{(i,j)\in\Sc}x^{i}y^{j}$ invariant. The question Q2 was originally answered in the following  two theorems.
\begin{thm}[see \cite{BMM2010}] \label{22_BMM}
For the 16 walks with a group of order 4 and for the 3 walks in figure~\ref{Fig_hol}, the 
formal trivariate series~\eqref{eq:tri} is holonomic non-algebraic. For the 
3 walks on the left in figure~\ref{Fig_alg}, the trivariate series~\eqref{eq:tri} is algebraic. 
\end{thm}
\begin{thm}[see \cite{BK2010}] \label{1_BK}
For the so-called \emph{Gessel's walk} on the right in figure~\ref{Fig_alg}, 
the formal trivariate series~\eqref{eq:tri} is algebraic. 
\end{thm}

Proving Theorem~\ref{22_BMM} requires skillful algebraic manipulations together with the calculation of adequate \emph{orbit} \index{orbit} and \emph{half-orbit} sums.
As for Theorem \ref{1_BK}, it has been mainly obtained by the powerful computer algebra system \emph{Magma}, which allows dense calculations to be carried out.
In \cite{FR2010}, a direct proof of these  theorems have been proposed by application of general results given in \cite{FIM2}, together with the fact that  $\pi(x)$ in equation~\eqref{eq:eqfonc} is always holonomic.  
\subsection{Explicit solutions and asymptotics (see \cite{FR2012})} Along the lines sketched in the preceding sections, it is possible to define a BVP for, say,  $F(x,0,z)$. Here, all the objects coming in the formulas depend on $z$, which merely acts as a complex parameter and  appears either as a subscript or an argument. 

The following formula is direct, since here the BVP is of Dirichlet Carleman type, due to the simple form of the coefficients $c(x)$ and $\widetilde{c}(y)$ in \eqref{eq:cgf}.

\begin{prop}
For $x\in \mathscr{G}(\M_z)$, 
\begin{equation} \label{eq:F(x,0,z)}
c(x)F(x,0,z)-c(0)F(0,0,z)=\frac{1}{2\pi i z}\int_{\M_z} tY_0(t,z)\frac{w'(t,z)}{w(t,z)-w(x,z)}\,\textnormal{d}t,
\end{equation}
where $w(x,z)$ is the gluing function for the domain $\mathscr{G}(\M_z)$ 
in the $\CC_x$-plane.
Of course, a similar expression could be written for $F(0,y,z)$.
\end{prop}
\subsection{On the singularities of the generating functions}
By symmetry and classical arguments,  we note  that only real singularities  of $F(0,0,z)$, $F(1,0,z)$, $F(0,1,z)$ and $F(1,1,z)$ with respect to $z$ will play a role in the asymptotics. From the expression \eqref{eq:F(x,0,z)}, the main origin of all possible singularities can be explained. We simply quote the main result (see \cite{FR2012}).

\begin{prop}
\label{prop:singularity-F(0,0,z)}
The smallest positive singularity of $F(0,0,z)$ is
\begin{equation}
\label{def_z_g}
     z_g = \inf \{ z>0 : y_2(z) = y_3(z)\}.
\end{equation}
\end{prop}

\begin{rem}
\label{rem:singularity-F(0,0,z)}
We chose to denote the singularity above by $z_g$, as one alternative definition  could be the following: the smallest positive value of $z$ for which the genus of the algebraic curve $\{(x,y)\in\mathbb C^2: K(x,y,z)=0\}$ \emph{switches from $1$ to $0$}.  In \cite{FR2012}, five equivalent characterizations of $z_g$ are proposed. 
\end{rem}

\subsection{The simple random walk}
\label{sec:simple-walk}
 For the \emph{simple walk} [i.e.  $(i,j)\in\Sc $ if and only if  $ij=0$], formulas are pleasant, because then the curve $\M_z$ is a \emph{circle}.
\begin{prop}
For the simple walk,
 \begin{eqnarray*}
     F(0,0,z) & =  & \frac{1}{\pi} \int_{-1}^1 \frac{  1-2uz-\sqrt{(1-2uz)^2-4z^2}}{z^2}
     \sqrt{1- u^2}\,\textnormal{d}u,\\
     F(1,0,z) & = & \frac{1}{2\pi}\int_{-1}^{1}\frac{1 - 2uz -\sqrt {(1-2uz)^2 - 4z^2}}{z^2}\sqrt{\frac{1+u}{1-u}}\, \textnormal{d}u.
\end{eqnarray*}
\end{prop}
Note that $F(0,0,z)$ counts the number of excursions starting from~$(0,0)$ and returning to~$(0,0)$, while $F(1,0,z)$ counts the number of walks starting from~$(0,0)$ and ending at the horizontal axis. By symmetry, $F(0,1,z) = F(1,0,z)$.
The following asymptotics holds.
\begin{prop}
\begin{equation*}
     f(0,0,2n)\sim \frac{4}{\pi}\frac{16^{n}}{n^3},\quad 
\sum_{i\geq 0}f(i,0,n)\sim \frac{8}{\pi}\frac{4^{n}}{n^2}, \quad
     \sum_{i,j\geq 0}f(i,j,n)\sim \frac{4}{\pi}\frac{4^{n}}{n}.
\end{equation*}
\end{prop}

\section{About generalizations}\label{sec:misc}
Some examples presented in this survey already contain some generalizations.  
There are essentially three main possible extensions.  First, for finite jumps of arbitrary size. Second, when the \emph{maximal space homogeneity} condition does not hold. 
Third, for random walks in $\Zb_+^n, n\ge3$.  
The reader will observe that these classes of problems are mathematically not disjoint. 
\subsection{Arbitrary Finite Jumps} 
Undoubtedly the first step toward a generalization, in the case of jumps bounded in modulus by a finite number $n$, is the analytic continuation process, which is crucial in most of the problems, including asymptotics.  Here there are $2n$
unknown functions, $\pi_i(s),\widetilde{\pi} _i(s)$, which must be
analytic in the connected domain $\mathcal{E}\subset\S$,
$$ \mathcal{E} = \{ |x(s)| <1,\,|y(s)|<1 \}.$$
 Then a functional equation can be obtained, on a Riemann surface $\S$ of
arbitrary genus, which has the form \index{genus ! arbitrary}
\begin{equation}
\sum_{i=1}^{n}\biggl(q_i(s)\pi_i(x(s))+\widetilde{q}_i(s)
\widetilde{\pi}_i(y(s))\biggr)+q_{0}(s)=0
\end{equation}
where $q_i(s),\widetilde{q}_i(s)$ are meromorphic on $\S$. Several results about analytic continuation were proved in \cite{MALY73}.

In the recent preliminary study \cite{FR2015},  finding and classifying 
branch-points and their associated cuts in the complex plane appear to be  two crucial issues. Indeed, the genus of the surface $\S$ is larger than~$1$, thus implying to manipulate hyperelliptic curves. The ultimate goal would be to set a generalized BVP on a single curve for a vector of analytic functions\,: this remains a doable challenge.

\subsection{Space inhomogeneity}
 Here, each situation is  peculiar. For instance, even if some explicit cases can be solved by reduction to a single FE (see e.g.,~\cite{FAY82}), however, most of the time, it will be necessary to deal with systems of functional equations, like in Section~\ref{sec:JSQ}. Also, it is often possible to write a \emph{non-N{\oe}therian} BVP (i.e. its index is \emph{not finite}) in some convenient regions.
\subsection{Larger dimensions} 
It is not necessary to insist on the usefulness of getting results for random walks in 
$\Zb_+^n, n\ge3$. Most of the questions are largely open. For a first step in this direction, see \cite{OVS}, where explicit integral formulas for the resolvent of the discrete Laplace operator in an orthant are obtained. At the moment, except for very speciaI cases, a global solution to the following problems seems  out of reach, even computationally\,: analytic continuation, index calculation, BVP for $n$ complex variables. 

The reason resides mainly in inductiveness properties: dimension~$n$ demands much finer properties for the related problems in dimension $n-1$, hence rendering the algebra almost untractable (even with the help of a computer programme\,!). But this is not too surprising, since, for instance, ergodicity conditions for random walks in $\Zb_+^n$ require finding invariant measures of walks in dimensions $\Zb_+^{n-1}$ !

The only tenuous hope might be to achieve a reduction to a vector BVP of a single variable on some hyperelliptic curve...

\end{document}